\documentclass[letter,scriptaddress,twocolumn, showkeys]{revtex4}
\usepackage{float} 
\usepackage{amsmath}
\usepackage{graphicx}
\usepackage{graphbox}
\usepackage{enumitem}
\usepackage{tikz}
\usepackage{tkz-euclide}
\usepackage{amssymb}
\usepackage{graphicx}
\usepackage[colorinlistoftodos]{todonotes}
\usepackage[utf8]{inputenc}
\usepackage{amsmath}
\usepackage[inline]{asymptote}
\usepackage{listings}
\usepackage{tabularx}
\usepackage{algorithm}
\usepackage{algpseudocode} 
\usepackage{fancyhdr}
\usepackage{booktabs}
\usepackage{color}
\usepackage{xcolor}
\usepackage{tikz}
\usetikzlibrary{arrows}
\usepackage{float}
\usepackage{comment}
\graphicspath{ {./} }
\usepackage[colorlinks=true, allcolors=blue]{hyperref}
\usepackage{setspace}
\usepackage{amsthm}
\newtheorem{thm}{Theorem}

\newtheorem{defn}[thm]{Definition}

\theoremstyle{remark}
\newtheorem{example}{Example}[]

\usepackage{multirow}
\usepackage{amscd,amssymb,latexsym,subfigure,verbatim,hyperref,epsfig,color}
 \usepackage{indentfirst}
      
\makeindex

\begin{document}
\title{CAACS: A {\it Carbon Aware} Ant Colony System}
\author{Marina Lin$^{a}$, 
 and Laura P. Schaposnik$^{\star, b}$ }
 \affiliation{($\star$) Corresponding author: schapos@uic.edu} 

\renewcommand{\thesection}{\Roman{section}}
\begin{abstract}
\par  In an era where sustainability is becoming increasingly crucial, we introduce a novel \textit{Carbon-Aware Ant Colony System (CAACS) Algorithm} that addresses the Generalized Traveling Salesman Problem (GTSP) while minimizing carbon emissions. This innovative approach harnesses the natural efficiency of ant colony pheromone trails to find optimal routes, balancing both environmental and economic objectives. By integrating sustainability into transportation models, the \textit{CAACS Algorithm} is a powerful tool for real-world applications, including network design, delivery route planning, and commercial aircraft logistics. Our algorithm's unique bi-objective optimization represents a significant advancement in sustainable transportation solutions.

\end{abstract}
\maketitle
\section{Introduction}
    
Sustainability is not just a buzzword but a critical goal that ensures our ability to safely co-exist on Earth for generations to come. Traditionally, sustainability encompasses three primary dimensions:
\begin{itemize}
    \item \textbf{Social} \cite{su9010068, VALLANCE2011342}
    \item \textbf{Economic} \cite{ELLIOTT2005263, Ikerd+2012}
    \item \textbf{Environmental} \cite{environmental1, environmental2}.
\end{itemize}
 As the scale and complexity of real-world problems grow, one of the biggest challenges is being able to balance multiple factors of sustainability. In this paper, we focus on the \textit{economic} and \textit{environmental} aspects. According to the United States Environmental Protection Agency \cite{EPA2024}, greenhouse gas (GHG) emissions from transportation account for 28\% of total U.S. greenhouse gas emissions, making it the largest contributor to U.S. GHG emissions. The largest sources of transportation include cars (at 48\%) and light and heavy-duty trucks (at 29\%).

One approach to mitigating GHG emissions is through optimizing transportation routes. The Traveling Salesman Problem (TSP) is a classical optimization problem that seeks the shortest possible route visiting a set of locations and returning to the origin point. First formulated in the 19th century \cite{1976Gt1, Korte2008}, TSP has since seen extensive algorithmic developments, achieving approximations within 1\% of the optimal solutions for problems involving millions of cities. The Generalized Traveling Salesman Problem (GTSP), an extension of TSP, is another NP-hard problem that has garnered significant research interest \cite{POP2024819}. Despite the extensive research on both GTSP and sustainability, to our knowledge, no existing GTSP algorithms explicitly account for sustainability considerations.
In this paper, we build on the concept of environmentally friendly TSP introduced in \cite{KAABACHI2017886} and extended in \cite{DAS2023101816} by developing a novel \textbf{Carbon-Aware Ant Colony System (CAACS)} Algorithm to solve the GTSP. This extension encapsulates a broader range of applications, offering a powerful tool for addressing sustainability in transportation models.

We begin this paper by recalling some background and introducing the Ant Colony System in Section \ref{sec:background}. We then define our novel carbon model in Section \ref{sec:carbon function}, inspired by \cite{MICHELI2018316}, which is a crucial part of our algorithm. The core of the paper, where we present a novel \textit{Carbon-Aware Ant Colony System Algorithm} to solve the Generalized Traveling Salesman Problem, appears in Section \ref{sec:CAACS}. In Section \ref{sec:scaling}, we discuss the effects of our novel emission factor, and  in Section \ref{sec:num_ants} the effect of the number of ants on the algorithm. We also conduct an empirical analysis on benchmark datasets of the GTSP in Section \ref{sec:empirical}. Finally, we analyze the time complexity in Section \ref{sec: time complexity}, leading to the following key results:

\begin{itemize}
    \item \textit{\textbf{Scaling Emission Factor:}} 
    \begin{itemize}
        \item Increasing the Scaling Emission Factor  $E(i, j)$ results in {\it shorter paths} up to a certain threshold;
        \item Higher values of $E(i, j)$ significantly {\it reduce carbon emissions}. 
    \end{itemize}
    \item \textit{\textbf{Ant Number:}} 
    \begin{itemize}
        \item A larger number of ants on the graph enhances solution quality.
        \item Increasing the number of ants reduces runtime up to a certain threshold. 
    \end{itemize}
\end{itemize}
 
The GTSP has a wide range of applications, including logistics, microchip design, scheduling, DNA sequencing \cite{DNA1, Dundar2019}, routing, material flow design, medical supply distribution, and image retrieval \cite{POP2024819}.  In Section \ref{sec: applications}, we demonstrate the novelty and utility of our algorithm by applying it to these diverse fields, showcasing its versatility and effectiveness.
\begin{itemize}
    \item \textit{\textbf{Network Design.}} 
    Our algorithm is used to create a more sustainable road network in the United States. We highlight specific characteristics that make the reduction of the GTSP to the TSP particularly effective for network design and similar applications.
    
    \item \textit{\textbf{Delivery Routes.}}    The algorithm optimizes delivery routes by identifying paths that minimize fuel consumption and reduce carbon emissions in logistics operations. In particular, we explore how our algorithm aligns with existing infrastructure and has the potential to accelerate drone-assisted parcel delivery services.
    
    \item \textit{\textbf{Commercial Aircraft Logistics.}} We discuss the application of our algorithm to commercial airplane logistics, including optimizing connecting flights and aircraft selection to improve efficiency and sustainability.\end{itemize}

\par We conclude this paper by discussing particularly notable aspects of the algorithm and results and outline future research directions in Section \ref{sec:conclusion}. 

\section{Background}
\label{sec:background}

In the following sections, we first present the Generalized Traveling Salesman Problem in Section \ref{subsec:GTSP}, previous approaches to the GTSP in Section \ref{subsec:previous approaches}, and finally discuss Ant Colony Optimization (ACO) in Section \ref{subsec: ACO}. Before we continue, we present Table \ref{table:symbol_table} to understand the notation used in the following sections.

\subsection{The Generalized Traveling Salesman Problem}
\label{subsec:GTSP}

The Traveling Salesman Problem has been a topic of great interest to mathematicians and computer scientists since the 19th century \cite{1976Gt1, Korte2008}. There are many variations of the TSP, but the one most relevant to this work is the Generalized Traveling Salesman Problem, which has many practical applications such as logistics, microchip design, scheduling problems, DNA sequencing \cite{DNA1, Dundar2019}, routing problems, material flow design problems, distribution of medical supplies, image retrieval, \cite{POP2024819}. The GTSP can more accurately model real-world problems than the TSP, especially situations with hierarchical structures (which can be modeled with clusters), creating more potential applications \cite{LIEN1993177}.

The General TSP is defined as follows. Consider an undirected graph $G = (V, E)$ whose vertices are divided into a given number of clusters, denoted by $C_1, C_2, \ldots, C_k$. The GTSP searches for the shortest Hamiltonian cycle (a path in an undirected or directed graph that visits each vertex exactly once) visiting a set of vertices such that exactly one vertex from each cluster is visited. The TSP is one specific instance of the GTSP where each cluster contains exactly one vertex, and it is an NP-hard problem (for details on previous solutions and challenges of the GTSP  see \cite{POP2024819}).
\begin{figure}[H]
    \centering 
    \includegraphics[width=0.4\textwidth]{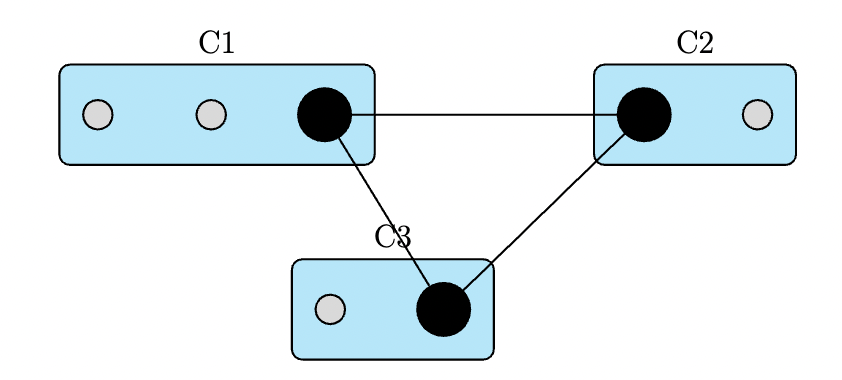}
    \caption{Example of a Valid Solution for a GTSP Problem with 3 clusters and 7 nodes.}
    \label{fig:gtsp_example}
\end{figure}
Currently, GLNS is the leading publicly available GTSP Solver \cite{SMITH20171}, based on the Lin-Kernighan heuristic, which is designed to be efficient and provide high-quality solutions for large instances. Additionally, OptimoRoute, a leading route optimization and scheduling software, uses algorithms that solve the Vehicle Routing Problem (VRP), which is an extension of the Traveling Salesman Problem (TSP). Both optimization software do not consider sustainability or offer alternative paths that are more environmentally friendly.

Before formulating the GTSP, we want to first mathematically formulate the original TSP problem, which can be either a symmetric TSP (sTSP), or an asymmetric TSP (aTSP). For the sTSP, let \( V = \{v_1,\ldots,v_n\} \) represent a set of cities, \( A = \{(r,s):r,s \in V\} \) be the edge set, and \( d_{rs} = d_{sr} \) be a cost measure associated with edge \( (r,s) \in A \). In this case cities \( v_i \in V \) are given by their coordinates \( (x_i,y_i) \) and \( d_{rs} \) is the Euclidean distance between \( r \) and \( s \) then we have an Euclidean TSP. The goal of the TSP problem is to find a minimal-length closed tour that visits each city once \cite{Matai10}.

\begin{defn}
 An aTSP occurs if \( d_{rs} \neq d_{sr} \) for at least one \((r,s)\).
\end{defn}

In order to introduce the formulation of the GTSP, consider the GTSP as a graph theoretic model. Let the graph \( G = (V, E) \) be a connected, undirected, and weighted graph, where each edge is associated with a non-negative cost. Without loss of generality, assume that G is a complete graph; however, if there is no edge between two nodes, assume there is infinite cost. The graph consists of a set of \( n \) vertices where \(V = 1, 2,   \ldots , n \) and a set of edges \( E = \{e_1, \ldots, e_m\} \) where $e_{i,j}=\{i,j\}\in E$ is the edge between nodes $i$ and $j$ and 
\[ E \subseteq \{ \{i, j\} \, | \, i, j \in V \text{ and } i \neq j \}. \]

To form the clusters,  the set of vertices \( V \) is divided into \( k \) disjoint nonempty subsets each of which forms {\it a cluster},  and which shall be denoted by \( C_1, C_2 \ldots, C_k \), so the following conditions hold:
    
\begin{enumerate}
    \item \(V =  C_1 \cup C_2 \cup \ldots \cup C_k \)
    \item \( C_i \cap C_p = \emptyset, \ \forall l, p \in \{1, \ldots, k\} \text{ with } l \neq p \).
\end{enumerate}

Associated with the model one has a cost function as \( c : E \to \mathbb{R}_+ \), which maps every edge \( e_{i, j} \in E \) of the graph to a positive number \( c_e = c(i, j) \in \mathbb{R}_+ \), called the cost of the edge \( e \). The associated cost matrix is denoted by \( W \). The edges are divided into intracluster edges that link two vertices from the same cluster and intercluster edges that link two vertices from different clusters.

\begin{defn}The GTSP is called Euclidean if the triangle inequality holds for any three vertices \(v_1, v_2, v_3 \in V\); otherwise, it is considered non-Euclidean. 
\end{defn}

\begin{defn} The GTSP is called symmetric if and only if the cost function \( c \) satisfies the equality \( c(i, j) = c(j, i) \) $\forall$ \( i, j \in V \).
\end{defn}

The GTSP is finding the minimum Hamiltonian tour \(H \) spanning a subset of nodes such that \(H \) contains exactly one node from each cluster \(C_i\) for \(i \in {1, \ldots, k}\). The L-GTSP problem involves \(H \) passing through each cluster at least once. In the specific instance of Euclidean problems, the optimal solution of the L-GTSP contains exactly one vertex from every cluster, reducing the problem to the GTSP with the same solution \cite{POP2024819}. For our project, we will be focusing on the first version of the GTSP defined,  where the cost matrix is symmetric, and exactly one vertex from every cluster is selected.\\

 \textbf{Mathematical Formulation of GTSP.} We shall follow the mathematical formulation of the standard GTSP in \cite{POP2024819} \cite{10.5555/2209505.2209517}.
Consider an undirected graph \( G = (V, E) \) and a subset of vertices \( S \subseteq V \). Let \[ E(S) = \{e_{i, j} \in E \,|\, i, j \in S\} \] represent the set of edges belonging to \( E \) with both vertices in \( S \). Define the cutest \( \delta(S) \) to be  \[\delta(S) = \{e_{i, j} \in E \,|\, i \in S, j \notin S\}, \] the set of edges in \(E \) where one vertex belongs to \(E \). To model the GTSP, let \(x_{ij}\) be a binary equal to $1$ if the salesman travels from node $i$ to node $j$ and $0$ otherwise, as seen below. The binary variable \(z_{i}\) is defined similarly: 
\begin{eqnarray}x_e = x_{ij} &=& 
  \begin{cases} 
   1 & \text{if the edge } e_{i, j} \in E \text{ is traveled on} \\
   0 & \text{otherwise}
  \end{cases}\nonumber
\\
z_i &=& \begin{cases}
1 & \text{if the vertex } i \text{ is visited}\\
0 & \text{otherwise}
\end{cases}\nonumber
\end{eqnarray}
 where \( i, j \in V \) and \( l, p \in \{1, \ldots, k\} \).

  One potential solution of the GTSP is a subgraph with no cycles and one vertex from each cluster, satisfying the degree constraints with a total of $k$ edges linking all the clusters. As a result, we can express the GTSP in terms of the binary variables:
\[
\begin{aligned}
& \min \sum_{e \in E} c(i,j) x_{i, j}
\end{aligned}
\]
subject to some constraints.     The first constraint ensures we are choosing exactly one vertex from each cluster:
    \[
    \sum_{i \in C_p} z_{i} = 1 \quad \forall p \in \{1, \ldots, k\}; \forall C_p \subseteq V.
    \]
  The second constraint guarantees that if a vertex is in the solution set, its degree must equal two:
    \[
    \sum_{\{i, j\} \in \delta(v)} x_{ij} = 2 z_{v} \quad \forall v \in V.
    \]
  The third constraint(s) are called the generalized subtour elimination constraints, which are used to exclude potential cycles. There is an exponential number of these constraints: 
    \[
    \sum_{\{i, j\} \in E(S)} x_{ij} \leq \sum_{k \in S \setminus \{i\}} z_k \quad \text{for } 2 \leq |S| \leq n - 2, \forall i \in S.
    \]
     The fourth constraint is the 0–1 integer constraint on the selection of edges within the solution:
    \[
    x_{e} \in \{0, 1\} \quad \forall e \in E.
    \]
   Finally, the fifth constraint is the 0–1 integer constraint on the selection of vertices within the solution: 
    \[
    z_{i} \in \{0, 1\} \quad \forall i \in V.
    \]

The formulation described above is called the \textit{generalized subtour elimination formulation}. The connectivity constraints may replace the subtour elimination constraints (3), as seen in \cite{POP2024819}, and the 0–1 linear formulation described is called the \textit{generalized cutset formulation}.

Since this formulation of the GTSP uses variables associated with the vertices and edges of graph $G$, they have an exponential number of constraints. A compact representation of the subtour elimination constraints uses a polynomial number of supplementary constraints and variables and is presented in \cite{10.5555/2209505.2209517} and \cite{pop2007}. 

Finally, the above GTSP reduces to create a complete tour $X_1, X_2, \ldots, X_N, X_1$ determined by the equation
\[
\text{min} \sum_{i=1}^{N-1} \sum_{j=1}^{N-1} c(X_i, X_j) + c(X_N, X_1), \quad i \neq j,
\]
 which is the objective function that minimizes the cost along with the criteria above. 

 \begin{table*}[ht]
\begin{tabularx}{\textwidth}{|l|X|l|X|}
\hline
\multicolumn{2}{|c|}{Parameters for Carbon Function \cite{MICHELI2018316}} & \multicolumn{2}{c|}{Parameters for Ant Colony Optimization} \\ \hline
\textbf{Symbol} & \textbf{Description} & \textbf{Symbol} & \textbf{Description} \\
$\xi$ & Fuel-to-air mass ratio & \(\tau_{ij}(0)\) & Pheromone on edge from node \(i\) to node \(j\) at the first iteration \\

$g$ & Gravitational constant & \( \mathcal{N}_{i}^{k}(t) \) & A set of feasible nodes to be visited by the \( k \)-th ant currently at node \( i \) \\

$\rho$ & Air density & \(\tau_{ij}(t)\) & Pheromone concentration on edge from node \(i\) to node \(j\) at iteration $t$ \\

$C_r$ & Coefficient of rolling resistance & $\rho$ & Evaporation rate\\

$\omega$ & Efficiency parameter for diesel engines & \(x^k(t)\) & The solution discovered by ant \(k\) \\

$\kappa$ & Heating value of typical diesel fuel &  \(f(x^k(t))\) & The quality of the solution found by ant \(k\) \\

$f$ & Vehicle speed & \(n_a\) & The number of ants \\

$\psi$ & Conversion factor & \( p_{ij}^{k}(t) \) & Transition Probability \\

$\varphi$ & Road angle & $\eta_{ij}(t)$ & The desirability of the move from $i$ to $j$, also known as \textit{visibility} \\

$\mu^k$ & Kerb-weight  & $d_{ij}(t)$ & Distance/Cost between edge $(i,j)$ \\

$c^k$ & Maximum payload (Capacity) &  $\Delta \tau_{ij}^k(t)$ & Pheromone Reinforcement\\

$k_e^k$ & Engine friction factor  & \(n_e\) & The number of elite ants\\

$N_e^k$ & Engine speed & $\hat{x}(t)$ & Current Best Solution at iteration $t$  \\

$V_e^k$ & Engine displacement & $f(\hat{x}(t))$ & Cost of Current Best Solution at iteration t \\

$C_d^k$ & Coefficient of aerodynamic drag & $x^+(t)$ & The best solution(s) giving the shortest path(s) \\

$A^k$ & Frontal surface area & $\hat{\hat{x}}(t)$ & The best path found from the first iteration to the current iteration $t$ of the algorithm \\

$\epsilon^k$ & Vehicle drive train efficiency & $f(x^+(t))$ & The cost(s) of the best solution(s) \\

$c(i, j)$ & Carbon Emission in kgCO$_2$e/litre & $r_0$ & Threshold set to balance exploration and exploitation \\

\( F_{i,j,k,p,t} \) & A payload in kg & $Q$ & Positive Scalar \\
\hline
\end{tabularx}
\caption{Symbol Table.}
\label{table:symbol_table}
\end{table*} 

\subsection{Approaches to the GTSP}
\label{subsec:previous approaches}

There have been several approaches to the GTSP that can be characterized into a few categories, namely 

\begin{itemize}
    \item Exact Algorithms \cite{henry1969record}
    \item Transformation Methods \cite{transformation}
    \item Reduction Algorithms \cite{reduction}
    \item Approximation Techniques \cite{approx}
    \item Heuristic Algorithms \cite{585892}
    \item Metaheuristic Algorithms \cite{pintea2017generalized}
\end{itemize}

A thorough comparison between various types of GTSP algorithms can be seen in \cite{POP2024819}. After considering the strengths and weaknesses of various algorithm structures, we found that sustainability can be most adaptable to metaheuristic algorithms which are generally used to solve complex optimization problems when a sub-optimal solution can be found \cite{DAS2023101816}.

One specific category of metaheuristics is Particle Swarm Optimization (PSO): similar to a swarm of animals, PSO consists of a swarm of agents who behave according to defined rules. For instance, any agent's movement could be determined by a mathematical formula that considers the positions of other agents as a way to find an optimal solution. Moreover,  PSO is able to search a large area in parallel and can model various animal behaviors \cite{Wang2018}.

\subsection{Ant Colony Optimization}
\label{subsec: ACO}

Ant Colony Optimization (ACO), as described in \cite{ACO}, is a model to find paths between nodes proposed in the late 20th century \cite{dorigo} inspired by the behavior of ant colonies. A key mechanism is that ants can leave behind pheromone trails, which other ants can sense and use to direct their behavior \cite{czaczkes2011synergy}. Using this information, an individual ant constructs potential solutions to a problem by adding components until a complete solution is created. Initially, \(\tau_{ij}(0)\), the pheromone on edge from node \(i\) to node \(j\) at the first iteration, is instantiated as a random small value to encourage exploration of different paths.

To select the next node to visit, ACO uses a \textbf{transition probability} known as the \textit{roulette wheel selection method}, which is defined as follows:  To select the next node \( j \in \mathcal{N}_i^k(t) \) from a set of feasible nodes to be visited by the \( k \)-th ant currently at node \( i \) for  \( \alpha > 0 \), the relative importance of the trail, the  {\it transition probability} is 
 
\begin{equation}
    p_{ij}^k(t) = 
    \begin{cases} 
        \frac{\tau_{ij}^\alpha(t)}{\sum\limits_{u \in \mathcal{N}_i^k(t)} \tau_{iu}^\alpha(t)} & \text{if } j \in \mathcal{N}_i^k(t) \\
        0 & \text{otherwise}
    \nonumber
    \end{cases}.
\end{equation}

 Note that if \( \mathcal{N}_i^k(t) = \emptyset \), the predecessor to node \( i \) is included in \( \mathcal{N}_i^k(t) \), which may cause loops. The loops are then removed when the destination has been reached.

\begin{example}\label{example1}
 To illustrate the evolution of path selection, consider this graph consisting of 4 nodes, and assume we start at $x_1$ and choose the closest node available. We are trying to find a valid TSP solution. The blue represents the node we are currently on, and the yellow represents nodes that have been selected. The dotted edges represent the edges that are in \( \mathcal{N}_i^k(t) \) at any given stage. 

\begin{figure}[H]
    \centering 
    \includegraphics[width=0.5\textwidth]{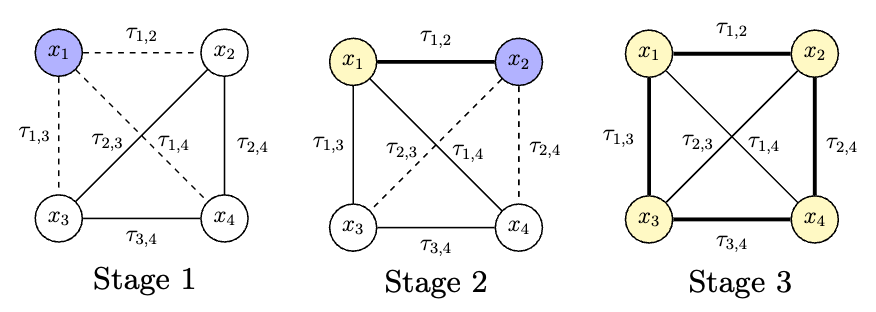}
    \caption{Example of path construction in a 4-node graph. Note that in reality, for ACO, the next node selected isn't always the closest, but the closer nodes have a higher probability of being selected. The thick edges represent the path selected. }
    \label{fig:example1}
\end{figure}
\end{example}

 \textbf{Update of Pheromone Intensity.} When updating the pheromone, the algorithm has two key steps: \textit{negative feedback}, where the pheromone is evaporated on all edges, and \textit{positive feedback}, where shorter paths are reinforced with more pheromone. 
 
During evaporation, for each edge \((i, j)\), pheromone intensity is reduced according to Eq.~\eqref{eq:eq1}:

\begin{equation}
    \tau_{ij}(t) \leftarrow (1 - \rho) \tau_{ij}(t), \label{eq:eq1}
\end{equation}

 \noindent where \(\rho \in (0, 1)\) is the evaporation rate. This allows ants to explore more and avoid premature convergence of the algorithm. To illustrate the importance of the evaporation process, we shall consider its effect on   Example \ref{example1}.

\medskip

\begin{example} Consider the setting of  Example \ref{example1}: we can see the effect of Eq.~\eqref{eq:eq1} on the 4 node pheromone graph as follows:

\begin{figure}[H]
    \centering 
    \includegraphics[width=0.4\textwidth]{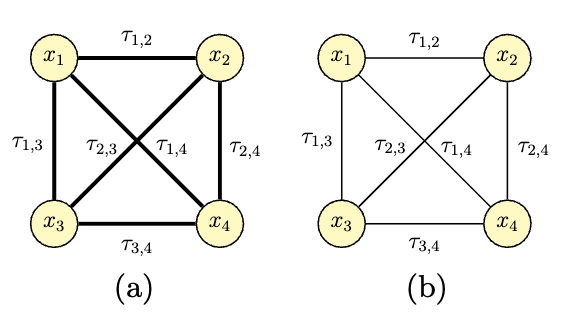}
    \caption{(a) Initial Pheromone buildup on a 4 node graph after multiple ants have found valid paths. (b) The removal of excess pheromones on the graph so the following ants can still explore new nodes.  }
    \label{fig:pheremone}
\end{figure}

After the above process, ants have constructed their paths from the source to the destination, all loops are removed, and the pheromone intensity on edge \((i, j)\) is adjusted to reinforce the shortest path: 

\begin{equation}
    \tau_{ij}(t + 1) = \tau_{ij}(t) + \sum_{k=1}^{n_a} \Delta \tau_{ij}^k(t), \label{eq:eq2}
\end{equation}
 where
  \[
  \Delta \tau_{ij}^k(t) =
  \begin{cases}
  \frac{Q}{f(x^k(t))} & \text{if edge } (i, j) \text{ occurs in path } x^k(t) \\
  0 & \text{otherwise}
  \end{cases}.
  \]

\begin{figure}[H]
    \centering 
    \includegraphics[width=0.4\textwidth]{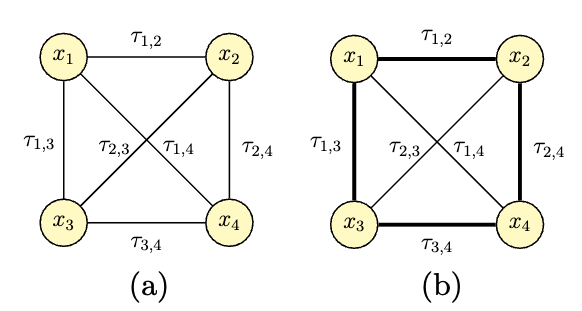}
    \caption{(a) Initial Pheromone distribution on a 4 node graph. (b) Pheromone Reinforcement on a 4 node graph of the shortest paths in Stage 3 of Figure \ref{fig:example1}.}
    \label{fig:pheremone2}
\end{figure}
\end{example}

Within this model, the shorter the path, the more additional pheromone. Components with greater pheromones will be more likely used in the following iterations of the algorithm. This model can be defined by an algorithm that is repeated many times by each ant: 

\begin{itemize}
    \item At every iteration, while the ant hasn't reached the node it started on, it continues to select the next node based on the transition probability above.
    \item After a path is found, the graph is updated with a new pheromone level depending on the path found.
\end{itemize} 

 The high-level process for Simple Ant Colony Optimization (SACO) can be seen in \textbf{Algorithm 1} and is expanded upon in \textbf{Algorithm 2}:

\begin{algorithm}[H]
\caption{The Ant Colony Optimization Metaheuristic \cite{ACO}}
\begin{algorithmic}[1]
\State Set parameters, initialize pheromone trails
\While{termination condition not met}
    \State Construct Ant Solutions
    \State Apply Local Search (optional)
    \State Update Pheromones
\EndWhile
\end{algorithmic}
\end{algorithm}

\begin{algorithm}[H]
\caption{Simple Ant Colony Optimization Algorithm \cite{ant_col_video1}}
\begin{algorithmic}[1]
\State Initialize $\tau_{ij}(0)$ to small random values; Let $t = 0$;
\State Place $n_a$ ants on the origin node;
\While{not STOP-CRITERIA}
    \For{each ant $k = 1, \dots, n_a$}
        \State $x^k(t) = 0$;
        \While{valid TSP solution not found}
            \State Select next node based on ACO $p_{ij}^k(t)$;
            \State Add $(i,j)$ to path $x^k(t)$;
        \EndWhile
        \State Remove all loops from $x^k(t)$;
        \State Calculate $f(x^k(t))$;
    \EndFor
    \For{each edge $(i,j)$ of the graph}
        \State Reduce the pheromone using Eq.~\eqref{eq:eq1};
    \EndFor
    \For{each edge $(i,j)$ of the graph}
    \State Update $\tau_{ij}$ using Eq.~\eqref{eq:eq2};
    \EndFor
    \State $t \gets t+1$;
\EndWhile
\end{algorithmic}
\end{algorithm}

 \noindent \textbf{Ant System.} The Ant System (AS) was developed based on SACO. In this model, there are multiple potential update strategies that are used to update pheromone intensity.  The main difference from the original ACO model is that heuristic information $n_{ij}$ is included in transition probability \( p_{ij}^k(t) \). 
 
 Here, $n_{ij} = \frac{1}{c_{ij}}$ is known as the \textit{visibility}. Additionally, a \emph{tabu list} to the set of feasible nodes, \( \mathcal{N}_i^k(t) \) may include only the immediate neighbors of the node \(i\) or may extend to all nodes not yet visited by ant \(k\) in order to prevent loops. Two different transition probabilities also characterize the Ant System where $\alpha \in [0, 1]$ is the relative importance of the trail, $\beta$ is the relative importance of the visibility, and $k \in \{1, \ldots, n_a\} $ is the current ant.

\medskip

 \textbf{Method 1:}
\[
p_{ij}^k(t) = 
\begin{cases} 
\dfrac{\tau_{ij}(t)^\alpha \cdot \eta_{ij}(t)^\beta}{\sum\limits_{u \in \mathcal{N}_i^k(t)} \tau_{iu}(t)^\alpha \cdot \eta_{iu}(t)^\beta} & \text{if } j \in \mathcal{N}_i^k(t) \\
0 & \text{otherwise}
\end{cases}.
\]

 \textbf{Method 2:}
\[
p_{ij}^k(t) = 
\begin{cases} 
\dfrac{\alpha \tau_{ij}(t) + (1 - \alpha) \eta_{ij}(t)}{\sum\limits_{u \in \mathcal{N}_i^k(t)} \alpha \tau_{iu}(t) + (1 - \alpha) \eta_{iu}(t)} & \text{if } j \in \mathcal{N}_i^k(t) \\
0 & \text{otherwise}
\end{cases}.
\]

In this model, pheromone evaporation occurs using Eq.~\eqref{eq:eq1}, and we update pheromone concentration according to Eq.~\eqref{eq:eq2} similar to SACO. However, $\Delta \tau_{ij}^k(t)$ changes depending on the Ant System algorithm used. There are three popular Ant System algorithms (see   \cite{ant_system} for further details) which we shall recall here and to which we shall come back later, where in  all three variations,  \( Q > 0 \) is a constant:

\begin{itemize}
  \item Ant-cycle AS: 
  \[
  \Delta \tau_{ij}^k(t) = 
  \begin{cases} 
    \dfrac{Q}{f(x^k(t))} & \text{if edge } (i,j) \text{ occurs in path } x^k(t) \\
    0 & \text{otherwise.}
  \end{cases}
  \]

  \item Ant-density AS: 
  \[
  \Delta \tau_{ij}^k(t) = 
  \begin{cases} 
    Q & \text{if edge } (i,j) \text{ occurs in path } x^k(t) \\
    0 & \text{otherwise.}
  \end{cases}
  \]

  \item Ant-quantity AS: 
  \[
  \Delta \tau_{ij}^k(t) = 
  \begin{cases} 
    \dfrac{Q}{d_{ij}(t)} & \text{if edge } (i,j) \text{ occurs in path } x^k(t) \\
    0 & \text{otherwise.}
  \end{cases}
  \]
\end{itemize}

Finally, an \textit{elitist strategy} is introduced in \cite{ant_system} to provide additional reinforcement to the best solution found so far by adding more pheromone to the paths walked by \textit{elite ants}, which is the ant that found the current best solution ($\hat{x}(t)$) at some iteration $t$. This cost of this solution is $$f(\hat{x}(t)) = \min\limits_{k=1,\ldots,n_a} f(x^k(t)).$$ This ensures that the optimal path's pheromone levels increase more than other solutions. The additional pheromone contribution per elite ant can be written as 
 \[\Delta \tau_{ij}^{e}(t) = 
  \begin{cases} 
    \dfrac{Q}{f(\hat{x}(t))} & \text{if edge } (i,j) \in \hat{x}(t) \\
    0 & \text{otherwise}
  \end{cases}\]

 which can then be added to Eq.~\eqref{eq:eq2} to get:

\begin{equation}
    \tau_{ij}(t + 1) = \tau_{ij}(t) + \sum_{k=1}^{n_a} \Delta \tau_{ij}^k(t) + n_e \Delta \tau_{ij}^{e}(t).\label{eq:eq3}
\end{equation}

It is important to note that the algorithm in \cite{DAS2023101816}, one of the first sustainable algorithms applied to the Traveling Salesman Problem, is based on the AS structure. Additional examples of ant algorithms for TSPs can be seen in \cite{JUNMAN2012319, app132111817}. 

 The algorithm for Ant System Optimization can be described as follows:

\begin{algorithm}[H]
\caption{Simple Ant System Algorithm \cite{ant_system}}
\begin{algorithmic}
\State Initialize $\tau_{ij}(0)$ to small random values, Let $t = 0$;
\State Set Parameters $\alpha, \beta, \rho, Q, \eta_{ij}(t)$;
\State Place $n_a$ ants on the origin node;
\While{not STOP-CRITERIA}
    \For{each ant $k = 1, \dots, n_a$}
        \State $x^k(t) = \varnothing$;
        \While{valid TSP solution not found}
            \State Select next node based on AS $p_{ij}^k(t)$;
            \State Add $(i,j)$ to path $x^k(t)$;
        \EndWhile
        \If{\textit{tabu list} is not used}
            \State Remove all loops from $x^k(t)$;
        \EndIf
        \State Calculate $f(x^k(t))$;
    \EndFor
    \For{each edge $(i,j)$ of the graph}
        \State Reduce the pheromone using Eq.~\eqref{eq:eq1};
    \EndFor
    \For{each edge $(i,j)$ of the graph}
        \If{elitist strategy is not used}
            \State Update using Eq.~\eqref{eq:eq2};
        \Else
            \State Update using Eq.~\eqref{eq:eq3};
        \EndIf
    \EndFor
    \State $t \leftarrow t+1$;
\EndWhile
\end{algorithmic}
\end{algorithm}

\smallbreak

 \noindent \textbf{Ant Colony System.} To better understand our proposed \textit{CAACS} Algorithm which we shall describe below, we introduce the Ant Colony System Algorithm (ACS) \cite{585892}, which has previously been proposed to solve the TSP and found that ACS outperformed AS in terms of solution quality and computational time. Local search methods have been implemented with ant colony optimization (ACO) to create hybrid algorithms that further improve the accuracy of the solutions, as seen in \cite{Karapetyan2012AnEH}. Later on, a Reinforcing ACS was developed to solve the GTSP \cite{Pintea2017}. A notable difference in the ACS compared to the Ant System is in the \textbf{Transition Probability} to select the next node.

\smallbreak

 \noindent \textbf{Transition Probability.} The \( k \)-th ($ k = 1, \ldots, n_a$) ant moves from node \( i \) to node \( j \) is according to:
\[
j = 
\begin{cases}
\underset{u \in \mathcal{N}_i^k(t)}{\mathrm{arg\,max}} \left\{ \tau_{iu}(t)^\alpha \eta_{iu}(t)^\beta \right\} & \text{if } r \leq r_0 \\
J & \text{otherwise}
\end{cases};
\]
 where \( J \in \mathcal{N}_i^k(t) \) is a node randomly selected according to the probability:
\[
p_{ij}^k(t) = 
\begin{cases}
\dfrac{\tau_{ij}(t)^\alpha \eta_{ij}(t)^\beta}{\sum\limits_{u \in \mathcal{N}_i^k(t)} \tau_{iu}(t)^\alpha \eta_{iu}(t)^\beta} & \text{if } J \in \mathcal{N}_i^k(t) \\
0 & \text{otherwise.}
\end{cases}
\]

ACS algorithm has two key components to its transition probability: \textit{exploitation}, where ants favor the best edge, and \textit{exploration}, where ants explore new unvisited nodes to find better paths. In the equation above, \( r \in [0,1] \) is a random number and when \( r \leq r_0 \) the algorithm exploits, and when \( r > r_0 \) the algorithm explores. The parameter \( r_0 \in [0,1] \) is used to balance exploration and exploitation.

The inclusion of \textit{visibility} ($\eta_{ij}(t)$) acts as additional heuristic information that guides ants to make decisions that are more informed rather than random. Furthermore, \textit{visibility} is used in conjunction with pheromone levels to balance the exploration of new paths and the exploitation of known good paths. While AS might favor more exploration, resulting in increased variability, ACS is likely to explore new paths when pheromone levels are low or equal and to follow stronger pheromone trails when \textit{visibility} is also high, thus reinforcing the better paths.

Another important feature of the ACS Algorithm is the inclusion of a \textbf{local and global update rule}. 
 {\it  The local update rule} is an extension of the evaporation rule in Eq.~\eqref{eq:eq1} that is applied to all edges where \( \rho_L \in (0, 1) \) a user-specified parameter and \( \tau_0 > 0 \) is a small constant:
\begin{equation}
    \tau_{ij}(t) \leftarrow (1 - \rho_L) \tau_{ij}(t) + \rho_L \tau_0. \label{eq:eq4}
\end{equation}
 {\it The global update rule}  is an extension of Eq.~\eqref{eq:eq2} to provide additional reinforcement of pheromone concentrations on the edges of the \textbf{best path} where $\rho_G \in (0, 1)$ is a user-specified parameter. It is important to note that at each iteration, $x^+(t)$ is equal to $\hat{x}(t)$ and is known as the \textit{iteration-best strategy}, but from the first iteration to the current iteration $t$ of the algorithm $x^+(t)$ is equal to $\hat{\hat{x}}(t)$ which is also known as the \textit{global-best strategy}.

\begin{equation}
    \tau_{ij}(t + 1) = (1 - \rho_G) \tau_{ij}(t) + \rho_G \Delta \tau_{ij}(t), \label{eq:eq5} 
\end{equation}
where
\[
\Delta \tau_{ij}(t) = 
\begin{cases} 
\frac{1}{f(x^+(t))} & \text{if } (i,j) \in x^+(t) \\
0 & \text{otherwise.}
\end{cases}.
\]
 
By including \textit{visibility}, ACO algorithms incorporate local problem-specific knowledge, which can help make more efficient and effective search decisions. Here, we present a standard ACS Algorithm implementation for the Traveling Salesman Problem. 

\begin{algorithm}[H]
\caption{Ant Colony System Algorithm \cite{585892}}
\begin{algorithmic}
\State Initialise $\tau_{ij}(0)$ to small random values, Let $t = 0$;
\State Set Parameters $\alpha, \beta, \rho_L, \rho_G, \tau_0, \eta_{ij}(t)$;
\State Place $n_a$ ants on the origin node, 
\State Initialize the best solution: $x^+(t) = \emptyset, f(x^+(t)) = 0$
\While{not STOP-CRITERIA}
    \For{each ant $k = 1, \dots, n_a$}
        \State $x^k(t) = \emptyset$;
        \While{valid TSP solution not found}
            \State Select next node based on ACO $p_{ij}^k(t)$;
            \State Add $(i,j)$ to path $x^k(t)$;
        \EndWhile
        \If{\textit{tabu list} is not used}
            \State Remove all loops from $x^k(t)$;
        \EndIf
        \State Calculate $f(x^k(t))$;
    \EndFor
    \For{each edge $(i,j)$ of the graph}
        \State Apply local update rule: Eq.~\eqref{eq:eq4};
    \EndFor
    \State Update global best solution: $x^+(t)$ and $f(x^+(t))$;
    \For{each edge $(i,j)$ in $x^+(t)$}
        \State Update global update rule: Eq.~\eqref{eq:eq5};
    \EndFor
    \State $t \gets t + 1$;
\EndWhile
\end{algorithmic}
\end{algorithm}

\section{Carbon Emission Function}
\label{sec:carbon function}

Since we are interested in the applications of our novel \textit{CAACS Algorithm}, we shall consider here carbon emissions from transportation. The International Energy Agency (IEA) has identified transportation as one of the largest contributors to carbon emission at 23\% \cite{su152115457}.  These emissions are dependent on a variety of factors such as the surface condition of the traveled route, conveyance type, load, speed, the weight of the vehicles, and thus it is important to consider these factors in our model. To simulate the carbon emission in the GTSP, we considered an emission function inspired by \cite{MICHELI2018316}, which includes three different components of the fuel consumption function respectively: 

\begin{itemize}
    \item The engine module, expressed as \(\lambda y (d_{i,j}/f)\), which is linear in the travel time. 
    \item The speed module, expressed as \(\lambda \gamma^k \beta^k d_{i,j}f^2\), which is quadratic in speed.
    \item The weight module, expressed as \(\lambda \gamma^k s (\mu^{k} + F_{i,j,k,p,t})d_{i,j}\), which is independent of the vehicle speed.
\end{itemize}

As a result, carbon emissions can be represented as:

\begin{scriptsize}
\[
c(i, j) = \lambda u \left( y \left( \frac{d_{ij}}{f} \right) + \gamma^ k \beta^k d_{ij} f^2 + \gamma^k s \left( \mu^k + F_{i,j,k,p,t} \right) d_{ij} \right)
\]
\end{scriptsize}

\noindent given a set of decision variables, general parameters, and vehicle type-dependent features. Additionally, the factor $u$ represents the fuel-dependent conversion factor, expressed in kgCO$_2$e/litre. In our model, $u$ equals 2.63 kgCO$_2$e/litre. Here, kgCO$_2$e represents kilograms of carbon dioxide equivalent, and our carbon emission is measured in kg of CO$_2$ emitted. In this model, we consider three vehicle types: Low-Diesel Vehicles (LDV), Medium Diesel Vehicles (MDV), and High Diesel Vehicles (HDV). The remaining constant values for each type of vehicle can be found in \cite{MICHELI2018316}. 

These three vehicle types made up 77\% of all global transport-related greenhouse gas emissions (GAE) in 2018 and are projected to make up 65\% of the transport-related GAE in 2060 \cite{ICCT2018}. Moreover, LDVs, including passenger cars and light trucks, constitute the majority of the vehicle fleet in the U.S. and account for over 90\% of all vehicles on the road. MDVs are less common than LDVs and makeup about 5-7\% of the total vehicle fleet. Heavy-duty vehicles (HDVs), including large trucks and buses, comprise about 3-5\% of the total vehicle fleet \cite{EDF2021}. Therefore, in the present paper, we shall consider only  LDVs, since they represent the majority of vehicles, but our code can be easily adapted to other vehicle types. 

\smallbreak

\noindent The variable parameters that we shall consider in our model will be 
\begin{itemize}
    \item \( f \) (m/s): a vehicle speed \( f \) (m/s).
    \item \( d_{i,j} \) (m): a traveled distance. 
    \item \( F_{i,j,k,p,t} \) (kg): a payload. This is representative of the correlation that heavier vehicles consume more fuel and produce more carbon emissions.
\end{itemize}

Since different types of highways and roads in the US have different speed limits and numerous paths between cities, it will be important to randomly select a speed between 11 m/s and 38 m/s for traveling between any two cities. This is supported by a study in 2015 that found that free-flow traffic on freeways averaged 31.46 m/s (70.4 mph), while major arterials averaged 25.20 m/s (56.4 mph), and minor arterials averaged 22.21 m/s (49.7 mph),  reflecting the diverse speed limits across different types of roads \cite{NHTSA2015}. The FHWA uses the 85th percentile speed to set speed limits, representing the speed at or below which 85\% of drivers travel under free-flow conditions. This method ensures that speed limits are safe and reasonable for most drivers. Speed limits on highways typically range from 24.59 m/s (55 mph) to 35.76 m/s (80 mph), depending on the state and type of roadway \cite{FHWA2013}. We shall also use realistic vehicle payloads, depending on the type of vehicle \cite{neighborhoodroadside2024, cleantechnol3020028}.
From the above analysis, our  general parameters shall be defined to be: 
\[ \lambda = \frac{\xi}{(k \psi)}       ~{\rm~ and~  }~  s = \tau + g \sin \phi + g C_r \cos \phi .\]

 The calculated vehicle type-dependent parameters are:
\begin{itemize}
    \item \( \gamma^k = \frac{1}{(1000 \omega \epsilon^k)} \).
    \item  \( \beta^k = 0.5 {C_d}^k \rho A^k \).
    \item \( y^k = {k_e}^k {N_e}^k {V_e}^k \).
    \item \( \mu^k \) : a kerb weight.
\end{itemize}

\section{CAACS: A Carbon Aware Ant Colony System Algorithm}
\label{sec:CAACS}

In what follows, we shall introduce our novel  \textit{Carbon-Aware Ant Colony System (CAACS) Algorithm}, which incorporates the carbon emission model above to solve the Generalized Traveling Salesman Problem. While this algorithm is an ant optimization algorithm inspired by those described in Section \ref{subsec:previous approaches}, there are some important differences we should note:

\begin{itemize}
    \item No previous algorithms include forms of sustainability to solve the GTSP. 
    \item  Ant Colony System has been known to produce higher-quality solutions compared to other previous Ant System-based algorithms that solely minimize cost \cite{DAS2023101816}.
    \item Our algorithm is easily adaptable to a variety of applications in multi-objective optimization. 
\end{itemize}

\begin{figure}[H] 
    \centering 
    \includegraphics[width=0.22\textwidth]{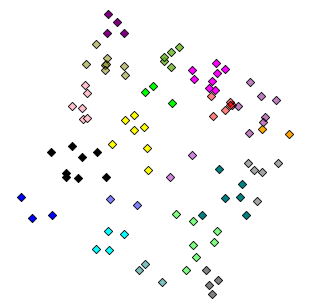} 
    \caption{Sample grid, with cities represented as 2x2
    diamonds and clusters as different colors. In this grid, there are 20 clusters and 100 nodes. Overlapping diamonds visually do not represent overlapping cities, but it is due to the scale of the figure and does not make a difference in the algorithm. } 
    \label{fig:sample_grid}
\end{figure}

 \noindent The search for a solution to a GTSP can be divided into two stages, which our algorithm shall take into consideration:
 \textbf{Stage I: Finding a Path} where a valid satisfies the GTSP and then \textbf{Stage II: Updating the Graph} where the graph updates the pheromone concentrations.  To illustrate our algorithm, we shall represent the cities, or nodes in the GTSP graph, as 2×2 diamonds on a grid, as shown in Figure \ref{fig:sample_grid}. The distance between each of the diamonds represents the relative distance between the cities. Each color represents a different cluster, so no two diamonds of the same color can be visited on the correct path. In this model, $N$ is the number of nodes, and $M$ is the number of clusters.

\subsection{Stage I: Finding a Path} In order to simultaneously minimize cost and emission, we introduce a novel emission scaling factor $E(i, j)$ below, which ranges between 0 (representing the least preferred paths with highest emissions) and 1 (most preferred paths with zero/minimal carbon emissions) 
\[ E(i, j) = A^{1 - \frac{C(i, j)}{C_{\text{max}}}} \]
where $C$ represents the carbon matrix and $A$ is a scaling factor. The \textit{CAACS Algorithm} will incorporate sustainability by integrating the carbon emission factor into the best path selection and reinforcement to select paths that minimize carbon emissions and cost for the Generalized Traveling Salesman Problem.  To do this, we shall start by incorporating the emission scaling factor into the transition probability to develop the \textbf{Carbon-Aware Transition Probability Rule}. This rule considers carbon emissions during path selection and consists of two main parts: \textit{Exploitation} and \textit{Exploration}. 

 \smallbreak

 \noindent \textbf{Exploitation.}
The \( k \)-th ant moving from node \( i \) to node \( j \) is according to
\[
j = 
\begin{cases}
\underset{u \in \mathcal{N}_i^k(t)}{\mathrm{arg\,max}} \left\{ \tau_{iu}(t)^\alpha \eta_{iu}(t)^\beta E_{ij} (t) ^\gamma \right\} & \text{if } r \leq r_0 \\
J & \text{otherwise.}
\end{cases}
\]
This occurs when $r \leq r_0$, which is a threshold that determines the balance of exploitation and exploration. The node $J$ is determined using the exploration below. 

 \smallbreak

 \noindent\textbf{Exploration.}
The node \( J \in \mathcal{N}_i^k(t) \),    for ant $k = 1, \ldots, n_a$, is randomly selected via the probability:
\[
p_{ij}^k(t) = 
\begin{cases}
\dfrac{\tau_{ij}(t)^\alpha \eta_{ij}(t)^\beta E_{ij}^\gamma (t)}{\sum\limits_{u \in \mathcal{N}_i^k(t)} \tau_{iu}(t)^\alpha \eta_{iu}(t)^\beta E_{ij}^\gamma (t)} & \text{if } J \in \mathcal{N}_i^k(t) \\
0 & \text{otherwise.}
\end{cases}
\]

\begin{figure}[H] 
    \centering 
    \includegraphics[width=0.22\textwidth]{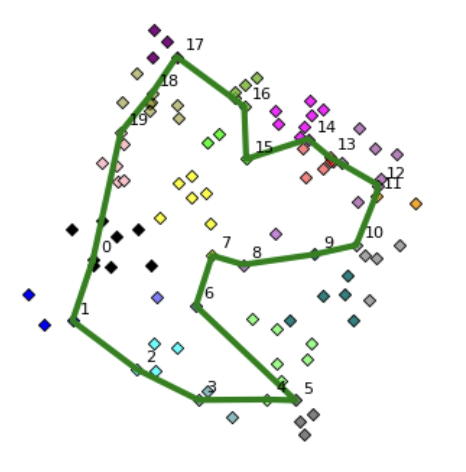} 
    \caption{An example of a valid solution found by an ant to the sample grid in Figure \ref{fig:sample_grid}.
    } 
    \label{fig:valid_sol}
\end{figure}

The transition probabilities favor edges with better environmental impact by increasing the pheromone for paths with less emission. As a result, ants have a higher probability of choosing a more environmentally beneficial path in their final solution. For the subsequent implementation of the algorithm, we want to balance exploration and exploitation so we set $r_0 = 0.5$.

\begin{example}

Given a GTSP graph with 5 clusters and 21 nodes, we visualize the impact of exploitation and exploration on segments of path creation.  

\begin{figure}[H]
    \centering 
    \includegraphics[width=0.36\textwidth]{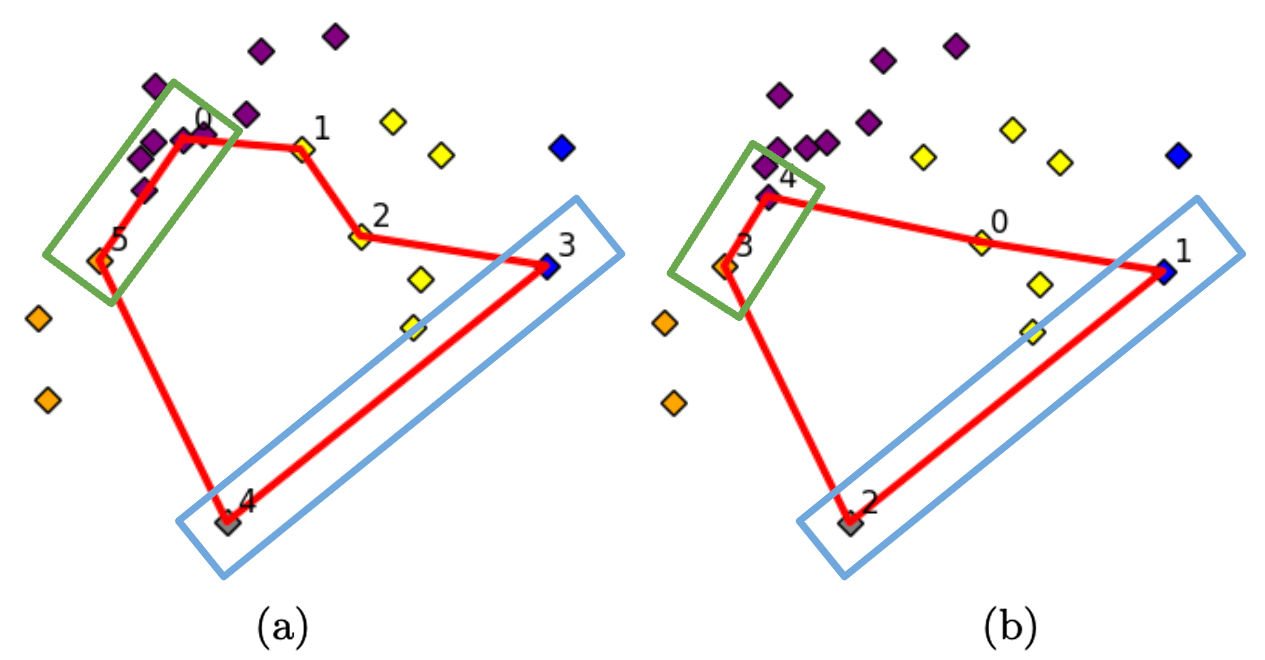}
    \caption{Example of the effects of (a) exploitation and (b) exploration on path selection. The green box represents one instance of the effect of exploration and discovery of a new path by a different ant at some iteration. The blue box represents the effect of exploitation and the selection of the same path between two nodes since it is optimal.}
    \label{fig:exploit}
\end{figure}

\end{example}

\subsection{Stage II: Updating the graph} In this portion of the algorithm, as it is often done in these types of algorithms \cite{engelbrecht2007computational}, updates to the pheromone levels in both the local and global update rules are part of the process to influence the path selection of the next ant (or ants in the next iteration). The local and global update rules incorporate sustainability by also adding the emission scaling factor $E(i, j)$ to previous models (e.g. see \cite{constantinou2010ant})

\noindent  \textbf{Carbon-Aware Local Update Rule.}
The local update rule is applied immediately after an ant traverses an edge and it is given by:
\begin{equation}
    \tau_{ij}(t) \leftarrow (1 - \rho_L) \tau_{ij}(t) + \rho_L \tau_0 E(i, j). \label{eq:eq6} 
\end{equation} It is designed to encourage exploration by reducing the pheromone level on the recently used edge, thereby making it less attractive for subsequent ants in the same iteration. 

Higher emission paths (lower $E(i, j)$) will contribute less to the pheromone update, making them less likely to be chosen in future iterations. Conversely, lower emission paths (higher $E(i, j)$) will have a stronger influence, promoting their selection.

\noindent \textbf{Carbon-Aware Global Update Rule.} The global update rule is applied after all ants have completed their tours and is given by:
\begin{equation}
\tau_{ij}(t + 1) = (1 - \rho_G) \tau_{ij}(t) + \rho_G \Delta \tau_{ij}(t) E(i, j)  \label{eq:eq7}
\end{equation}
where
\[
\Delta \tau_{ij}(t) = 
\begin{cases} 
\frac{1}{f(x^+(t))} & \text{if } (i,j) \in x^+(t) \\
0 & \text{otherwise.}
\end{cases}.
\]
It aims to reinforce the best solutions found so far, making them more attractive paths for future iterations. Once again, paths with lower emissions (higher $E(i, j)$) will receive a higher increase in pheromone, making the path more likely to be selected in the following iterations. 

\begin{figure*}[ht] 
    \centering 
    \hspace{-5.5in}
    \includegraphics[width=0.8\textwidth]{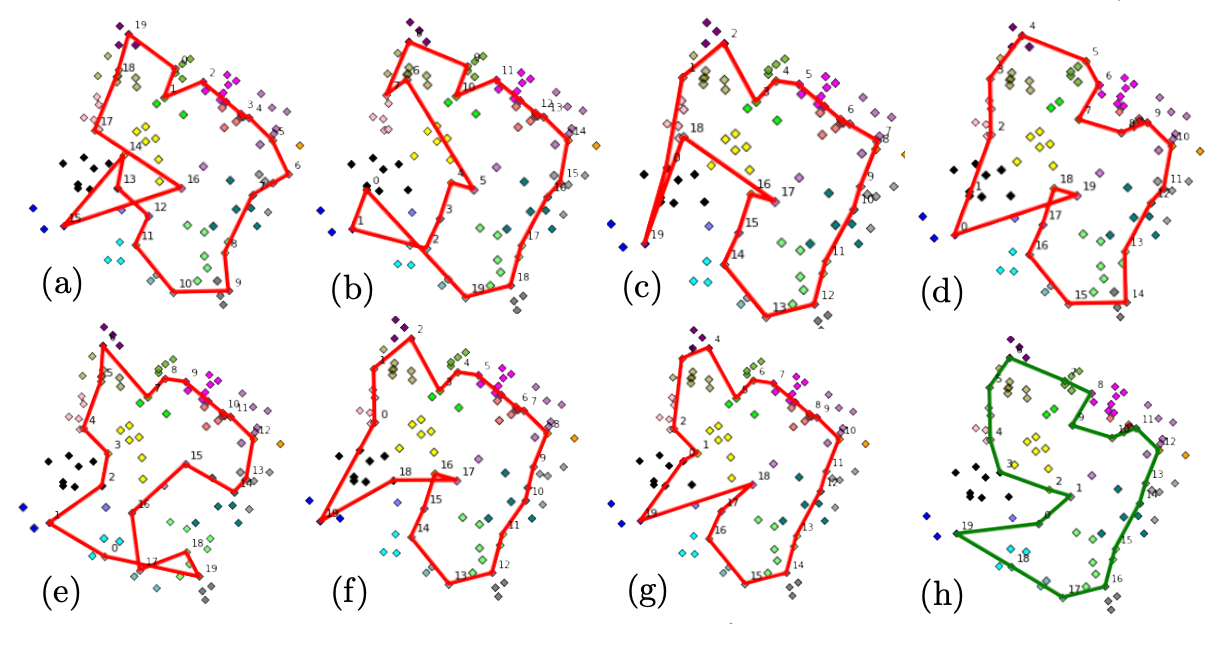} 
    \caption{Progression of the best path in the algorithm: (a) The first solution found by an ant. (b)-(g) Represent improved solutions found in subsequent iterations of the algorithm.  (h) the best solution found after all iterations are completed and a terminated criteria is met. In this case, the number of iterations is 100 and we are using the grid from Figure \ref{fig:sample_grid}.} 
    \label{fig:progression}
\end{figure*}

We shall also set a threshold on the maximum number of iterations in various trials that is always reached unless a terminating condition is specified in which case the algorithm ends. When updating a graph, segments of each image common to subsequent images represent the effects of {\bf exploitation} and new path segments represent the effects of {\bf exploration}. In what follows, we shall illustrate these updates by considering the progression of our algorithm applied to the network in  Figure \ref{fig:sample_grid}, and describe in Figure \ref{fig:progression} the progression. With the above definitions and analysis in mind, our  \textit{CAACS Algorithm} is defined as follows.

\begin{algorithm}[H]
\caption{Carbon Aware Ant Colony System Algorithm (CAACS)}
\begin{algorithmic}
\State Initialize $\tau_{ij}(0)$ to small random values, Let $t = 0$;
\State Set Parameters $\alpha, \beta, \rho_L, \rho_G, \tau_0, \eta_{ij}(t)$;
\State Initialize nodes, clusters, visited\_clusters \textit{tabu lists}
\State Place each of the $n_a$ ants on a random node
\State Initialize the best solution: $x^+(t) = \emptyset, f(x^+(t)) = 0$;
\While{not STOP-CRITERIA}
    \For{each ant $k = 1, \dots, n_a$}
        \State $x^k(t) = \emptyset$;
        \While{valid GTSP solution not found}
            \State Select next node based on CAACS $p_{ij}^k(t)$; 
            \State {// next node not in a visited\_clusters}
            \State visited\_clusters add the cluster of the next node
            \State Add $(i,j)$ to path $x^k(t)$;
        \EndWhile
        \State Calculate $f(x^k(t))$;
    \EndFor
    \For{each edge $(i,j)$ of the graph}
        \State Apply local update rule Eq.~\eqref{eq:eq6};
    \EndFor
    \State Update global best solution: $x^+(t)$ and $f(x^+(t))$;
    \For{each edge $(i,j)$ in $x^+(t)$}
        \State Update global update rule Eq.~\eqref{eq:eq7};
    \EndFor
    \State $t \gets t + 1$;
\EndWhile
\end{algorithmic}
\end{algorithm}

 \textbf{Parameters.} In this model, there are a few parameters to consider that are largely independent of the problem: $\beta, \rho_l, \rho_g$, and $\alpha$ introduced previously. In the present paper, all experiments are run with $\beta = 1$, $\rho_l = 0.99$, $\rho_g = \alpha = 0.1$. Other parameters in Table \ref{table:symbol_table} are specified before each trial.

\section{Scaling Emission Factor}
\label{sec:scaling}

In order to illustrate how different aspects of our model influence our algorithm, in what follows we shall consider how changing the emission scalar base $A$ may influence the results one may obtain through the \textit{Carbon-Aware ACS Algorithm} described in Section \ref{sec:CAACS}. 

We generated a new diverse dataset of 55 new GTSP instances with a variety of total nodes and clusters. For all the instances, the number of clusters $M$ is set following Eq.~\eqref{eq:eq8},

\begin{equation}
M = \frac{N}{5}. \label{eq:eq8}
\end{equation}

\noindent where the number of nodes $N$ ranges from 20 to 300. We generate 5 instances for each value of $N$, and we run 5 trials on each of the graphs given an emission scalar base $A$, for a total of 25 trials for every pair of $(M, N)$.
We measure the number of iterations that the algorithm takes to meet a terminating condition. This condition makes a discovered solution the final solution if it doesn't change after $\frac{N}{5}$ iterations or 1000 iterations are completed. We instantiated $30$ ants on the graph.

\begin{figure}[h]
    \centering 
    \hspace{-3.5in}
    \includegraphics[width=0.45\textwidth]{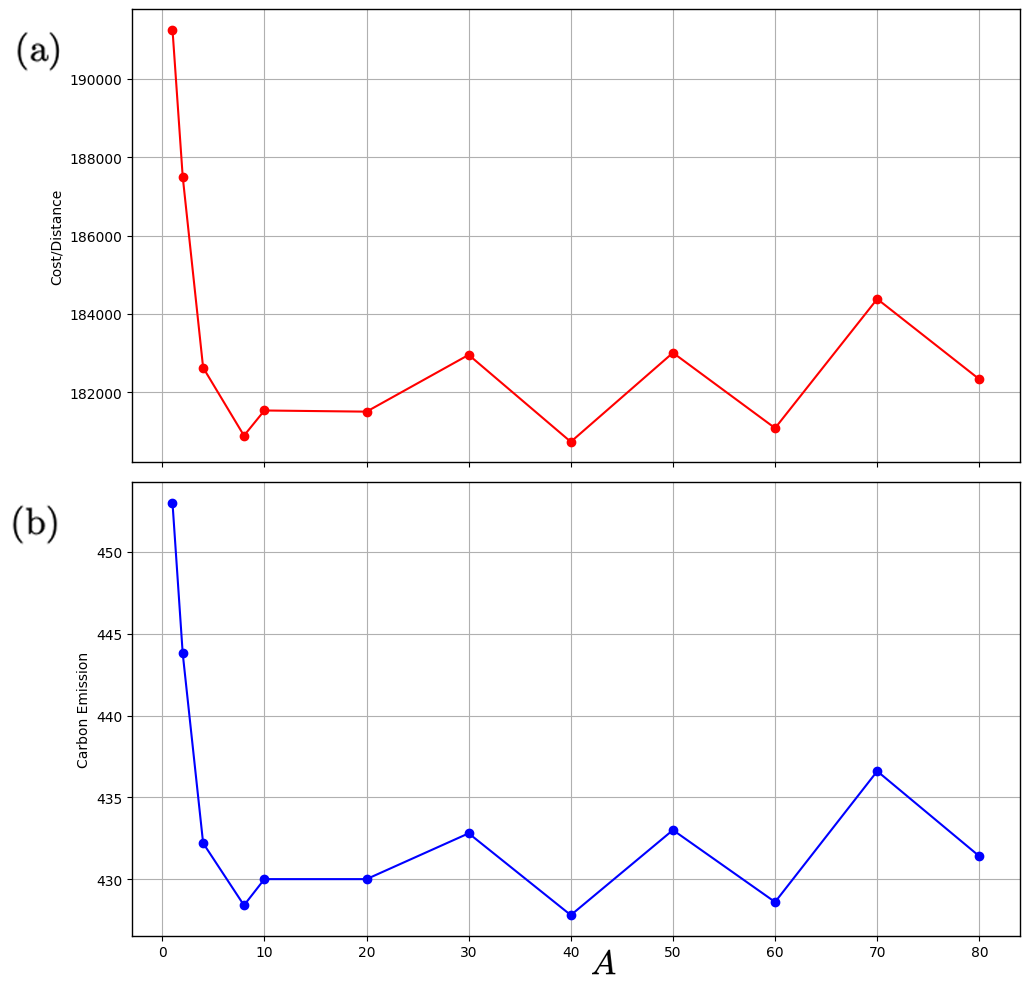}
    \caption{(a) Total cost and (b) Total carbon emissions varying with changes in $A$ for the dataset specified above, with  5 trials for each project instance, where we consider a fifth of the total cost and carbon emissions for all instances as a measure of solution quality.
    }
    \label{fig:combined_plots_random}
\end{figure}

As the value of $A$ increases, there is a decrease in carbon emissions, and this can be seen in Figure \ref{fig:combined_plots_random} below -- this is as expected, since an increase in $A$ increases the influence of the emission scalar factor in path selection and pheromone reinforcement. As a result, paths with lower carbon emissions are selected, decreasing the total carbon emission of the sum of all instances. One should note that setting A = 40 appears to be the most optimal in this case. 

Additionally,  as the value of $A$ increases, there is also a decrease in total cost, indicating the solution gets better. This is likely because the longer the path means increased carbon emissions, where the cost (or distance) is factored into the carbon function (however, it is not the only factor as described in Section \ref{sec:carbon function}). As $A$ increases, paths with lower carbon are favored, which also contributes to selecting paths with shorter costs. As a result, we see a similar trend in cost as we do in carbon.

\begin{figure}[h]
    \centering
    \hspace{-3.5in}
    \includegraphics[width=0.45\textwidth]{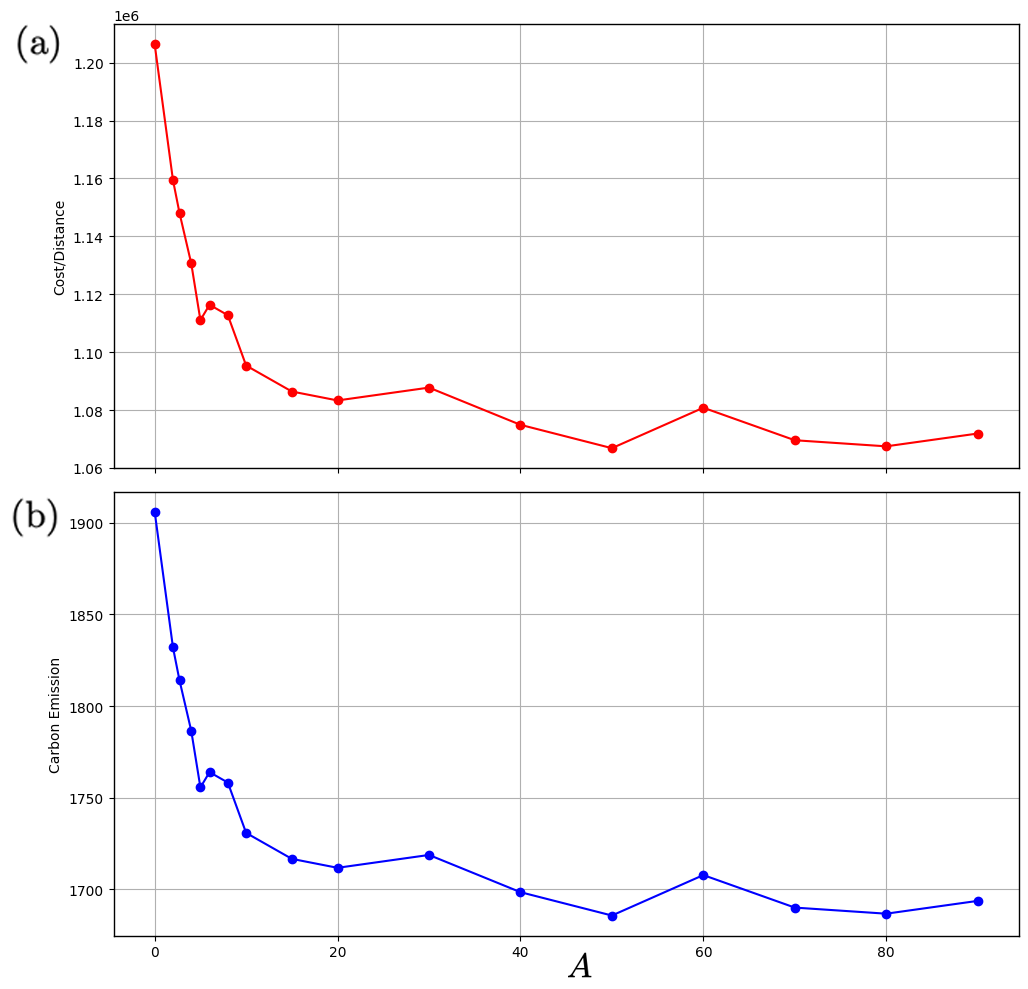}
    \caption{ (a) Total cost and (b) Total carbon emissions. $A = 0$ is calculated using $E(i, j) = 1$ which indicates the carbon emissions of the original GTSP.}
    \label{fig:combined_plots}
\end{figure}

To corroborate our findings further, we consider the benchmark GTSP Instances Library  \cite{gutin2009memetic}. We run 62 instances for every constant $A$. At each value of $A$, we took the sum of the cost and carbon emission for all instances as a measure of solution quality. During these trials, there were $30$ ants, and each instance was run for $100$ iterations. As shown in  Figure \ref{fig:combined_plots}, we see a similar trend to Figure \ref{fig:combined_plots_random} with a sharp improvement in cost and carbon within $A \leq 20$. Moreover, trials with $A = 50$ found the lowest carbon emission and cost (lowest data point between 40 and 60) which is within the range specified above. For larger values of $A$, there seems to be less fluctuation in total cost, most likely due to more realistic instances in a benchmark dataset as opposed to random generation of cost matrix, suggesting that  $ 40\leq A \leq 60$ would be approximately the best range to select $A$ from. 
 
\begin{figure}[h]
    \centering 
    \hspace{-2.5in}
    \includegraphics[width=0.35\textwidth]{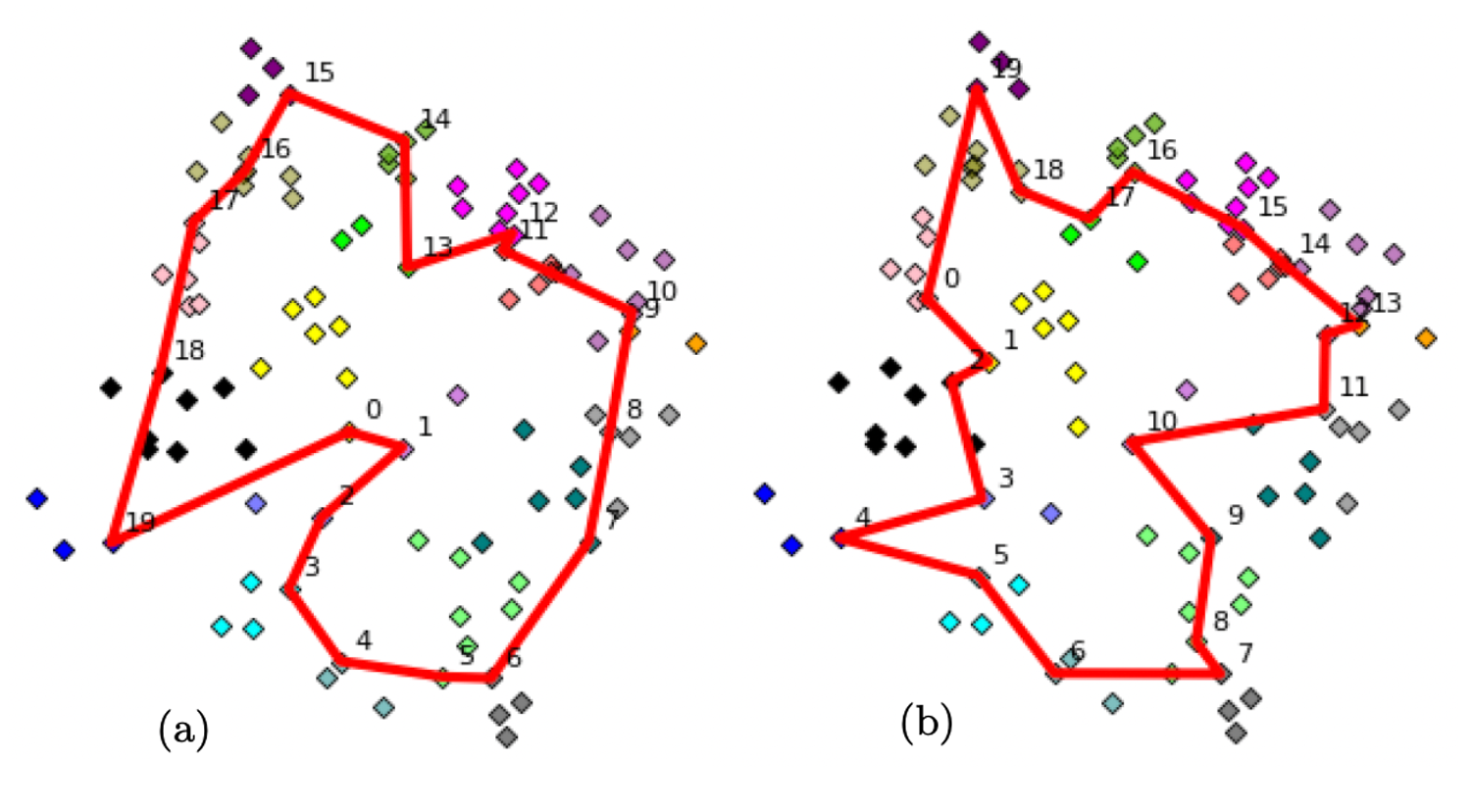}
    \caption{(a) Final solution  for $A = 10$ with 100 iterations with Carbon Emission = 9.198 kg CO$_2$ and Cost = 3927. (b) Final solution for $A = 50$ with 100 iterations, with Carbon Emission = 9.148 kg CO$_2$ and Cost = 3882. }
    \label{fig:a_comparison}
\end{figure}

Returning to the example considered in Figure \ref{fig:sample_grid}, one can see that 
  different $A$ results in a different final path generated, as shown in Figure \ref{fig:a_comparison}. In what follows, unless otherwise specified, we use $A = 50$ to achieve the optimal cost and carbon emissions. 

\section{The Effect of the Number of Ants}
\label{sec:num_ants}

We shall examine here the effects of the number of ants involved in the \textit{CAACS Algorithm}. We run 10 trials on 6 grids (varying between 100, 150, and 150 nodes and $\frac{N}{5}$ clusters) and every number of ants from 5 to 50 at intervals of 5. The six instances we use are from the GTSP Instances Library  \cite{gutin2009memetic}, and for each trial, we measure the cost discovered after 100 iterations. We use $r_0 = 0.5$ to balance \textit{exploration} and \textit{exploitation}. 

To simplify notation, we shall do our analysis by considering the triple $(M, N, id)$ given by

\begin{itemize}
    \item $M:=$  The number of clusters.
    \item $N:=$  The number of nodes.
    \item $id:=$  A unique identifier for the project instance. 
\end{itemize}

Given an example GTSP instance 20kroA100, we represent this as a tuple (20, 100, kroA). For project instances in this benchmark dataset that have a unique pair (M, N), we omit the third unique $id$ for simplicity. An example of this is 5ulysses22 and we denote this as (5, 22) since there is no other case in this set with 5 clusters and 22 nodes.

\subsection{Solution Quality}

One of the most important characteristics of our algorithm which changes with the number of ants is the quality of the solution. We measure this by considering the percentage error from the optimal solution provided in the GTSP Instances Library. We also compare the cost of the solution created by \textit{CAACS Algorithm} to the solution from our baseline GTSP model and observe the carbon improvement.

\begin{figure}[h]
    \centering 
    \hspace{-2.5in}
    \includegraphics[width=0.40\textwidth]{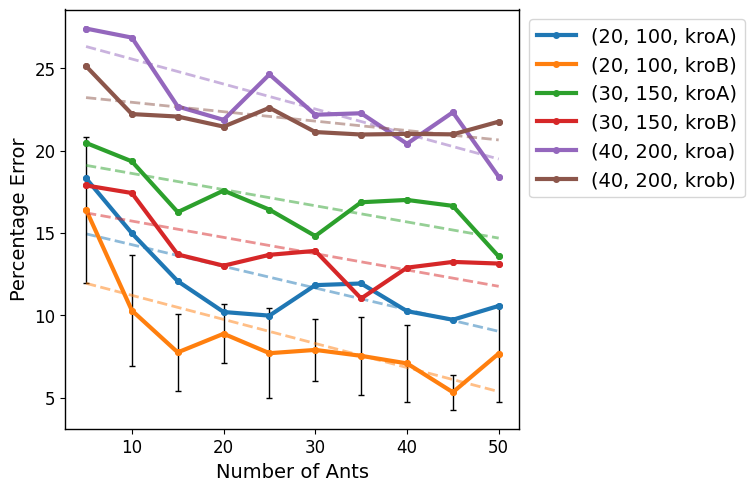}
    \caption{Percentage error in cost (quality of the solution) compared to the optimal solution versus the number of ants. Standard deviation showed for (20, 100, kroB). }
    \label{fig:num_ants_error}
\end{figure}
With more ants, the algorithm can explore a larger portion of the solution space which can increase the probability of finding a better solution since more paths are evaluated. There is also improved information sharing due to more pheromone trails laid down which can reinforce good solutions. As a result, more ants capitalize on previously discovered paths. These effects are depicted in Figure \ref{fig:num_ants_error} where one can see that overall, as the number of ants increases, the quality of the solution generally seems to improve.

It is important to note, however, that after 15 ants there is marginal improvement in the solution quality in some cases where solution quality worsened such as (40, 200, krob) from 20 to 25 and (30, 150, kroA) from 15 to 20 ants. This could be because the additional pheromone deposition could be over-saturated, making it difficult for following ants to differentiate between good and worse solutions. As a result, some ants would follow suboptimal paths, as seen in Figure \ref{fig:num_ants_error}.

To examine the correlations between the number of ants, cost, and carbon of our baseline GTSP ACS Algorithm and the \textit{CAACS Algorithm} we consider our algorithm through 10 trials on 6 grids and every number of ants from 5 to 50 for both models. Here, all instances of cost and carbon error are computed by calculating the percentage error between the cost/carbon respectively of each of the 10 trials using \textit{CAACS Algorithm} and the average cost/carbon of 10 trials of the baseline model. The negative values represent the cost or carbon decreased, which is essentially, the cost or carbon saved using our model. 

\begin{figure}[h]
    \centering 
    \hspace{-3.5in}
    \includegraphics[width=0.50\textwidth]{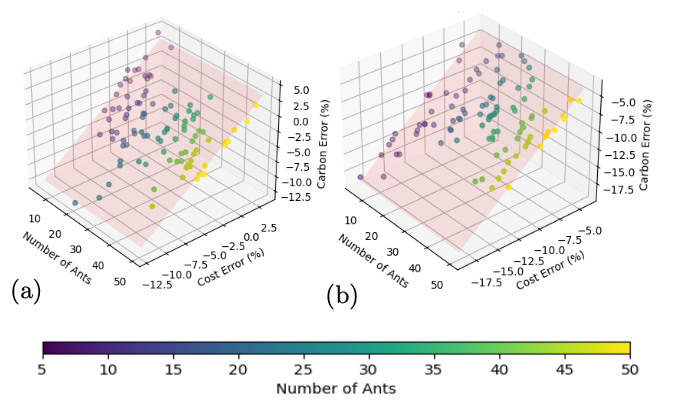}
    \caption{3D Scatter Plot Visualize Correlations between Number of Ants, Cost Error, and Carbon Error: (a) (20, kroA, 100) Instance. (b) (40, krob, 200) Instance. The red plane represents a best-fit plane for each 3D plot. Color Bar displays the color of points for a given number of ants. }
    \label{fig:3d_scatter}
\end{figure}

Our analysis is shown in Figure \ref{fig:3d_scatter}: we begin by considering the graph for a constant number of ants. We found in all instances, a nearly linear correlation between the cost error and carbon error.  Moreover, we see that for smaller instances (100 nodes), there is an even greater amount of carbon conserved than cost conserved since the slope is $\leq 1$. For larger instances, we find the carbon conserved is slightly less than the cost conserved for 150 nodes and 200 nodes.

Finally, we also find the slope in the $x$ directions from instances in (a)-(f) to be all relatively close to 0. Since the slope values are relatively close to zero indicates that there is a negligible correlation between the number of ants ($x$-axis) and the cost error ($y$-axis) for these instances.

\subsection{Run Time} 

\begin{figure}[H]
    \centering 
    \hspace*{-1.0cm}
    \includegraphics[width=0.45\textwidth]{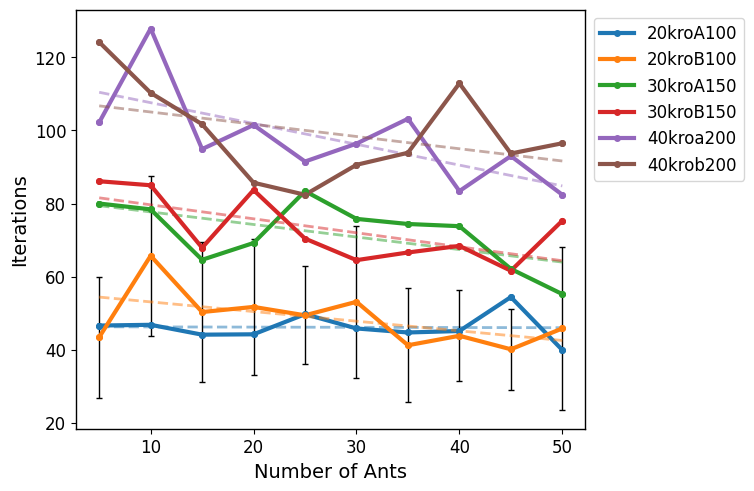}
    \caption{The number of iterations versus the number of ants. Standard deviation showed for (20, 100, kroB). }
    \label{fig:num_ants_iter}
\end{figure}

Finally, we shall consider the run time of our algorithm (the number of iterations the algorithm runs for) by setting termination criteria as if a solution found doesn't change after $\frac{N}{5}$ iterations, where  $N$ is the number of nodes. 
As the number of ants increases, one can see that the iterations required generally decrease then increase, particularly for medium/larger problem sizes (e.g. (30, 150, kroA) and (40, 20, krob)), as illustrated in Figure \ref{fig:num_ants_iter}. 
One may argue that increasing the number of ants improves the search efficiency of the algorithm for problems, but at a certain point, the potential for path redundancy of ants might reduce the overall effectiveness. We find that smaller problem sizes such as (20, 100, kroA) and (20, 100, kroB) typically have a stable number of iterations. However, for larger problem sizes ((40, 200, kroA) and (40, 200, krob), the relationship is less consistent, with some peaks and troughs. Additionally, high randomness in where the ants are initially placed can influence the overall number of iterations needed to converge.

From the above analysis, we see that using more ants allows one to explore paths and in general find shorter solutions, but can also become computationally more expensive past a certain threshold which we identify to be around 30 ants. As a result, we shall use 30 ants for the following trials to balance the quality of the solution discovered and the computation required. 

\section{Empirical Analysis}
\label{sec:empirical}

In what follows, we will conduct an empirical study of \textit{CAACS Algorithm} on the GTSP Instances Library  \cite{gutin2009memetic}, which is derived from TSPLIB instances. In the context of the GTSP, we denote $M$ as the number of clusters and $N$ as the number of nodes. There are 88 testing cases and the number of clusters ranges up to $M = 217$ and the number of nodes ranges up to $N = 1084.$

\begin{figure}[h]
    \hspace{-3in}
    \includegraphics[width=0.4\textwidth]{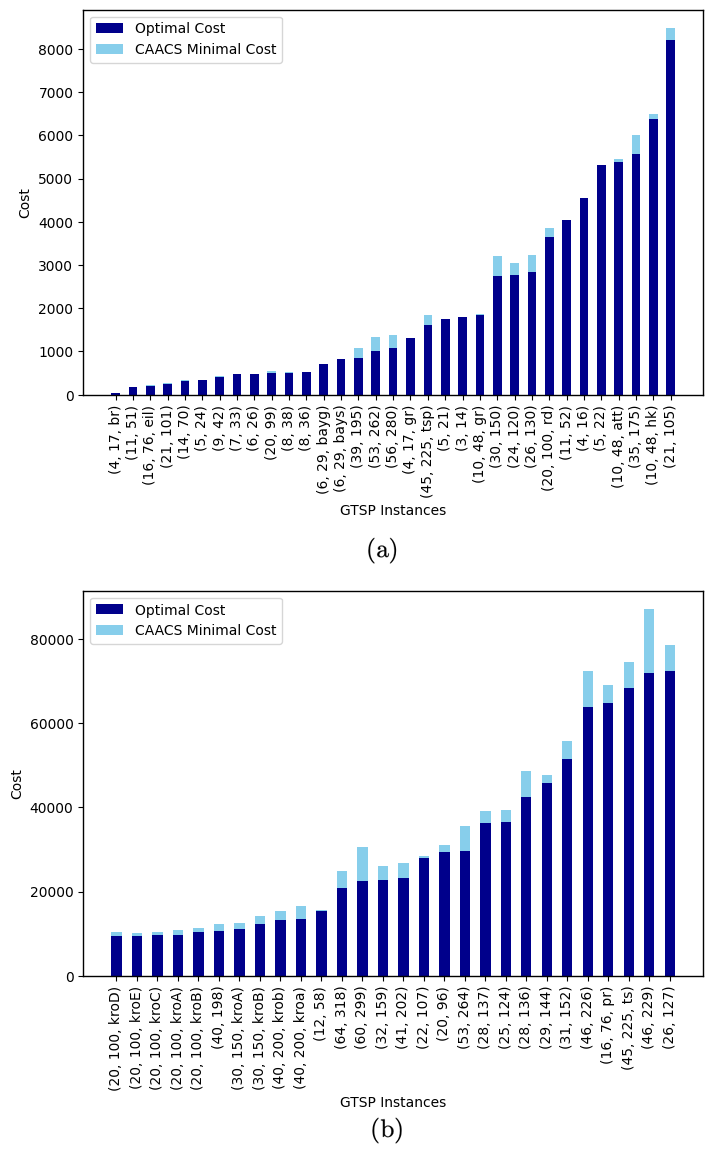}
    \caption{Comparison of optimal cost (in dark blue) and the \textit{CAACS} Minimal Cost (light blue). }
    \label{fig:exp_cost}
\end{figure}

For our experiment, we simulate 62 out of the 88 instances, up to $M = 64$ and $N = 300$. The selected 62 instances accurately represent the effects of the \textit{CAACS Algorithm} on GTSP on smaller instances and as the GTSP scales. Larger instances that are significantly more computationally expensive should scale according to the insights discussed below, allowing for the generalizability of our findings. We ran 10 trials on each instance and took the average cost of the 10 trials to compare with the optimal solution provided in \cite{gutin2009memetic}. In Figure \ref{fig:exp_cost} (a), we see that the optimal solutions are usually found consistently for small instances. As the optimal solution length increases, the error tends to increase, but not consistently, likely due to other factors like the number of nodes $N$ and the number of clusters $M$. Hence an increase in $N$ and $M$ would decrease solution quality, as instances such as (53, 262) and (56, 280) that have a similar optimal solution length to 4gr17. However, due to more nodes and clusters, there are increased constraints that need to be satisfied leading to more variability in the discovered solution.

Note that a portion of the light blue lies under the dark blue so \textit{CAACS} Minimal Cost can be seen as the sum of the dark blue and light blue. The light blue displayed would represent the error of the \textit{CAACS} in cost. The GTSP instances are ordered from least to greatest based on optimal cost. The instances are separated into two graphs to account for the wide range of costs. 

In Figure \ref{fig:exp_cost} (b) we see larger errors due to the more complex testing cases. We see that there are some outliers to the trend mentioned above: one example is that (60, 299) has fewer nodes and clusters than (64, 318) but there is a greater error. We attribute these other factors such as variables in path selection such as the influence of carbon emission and individual path distance. 

\begin{figure}[h]
    \centering 
    \hspace{-3in}
    \includegraphics[width=0.4\textwidth]{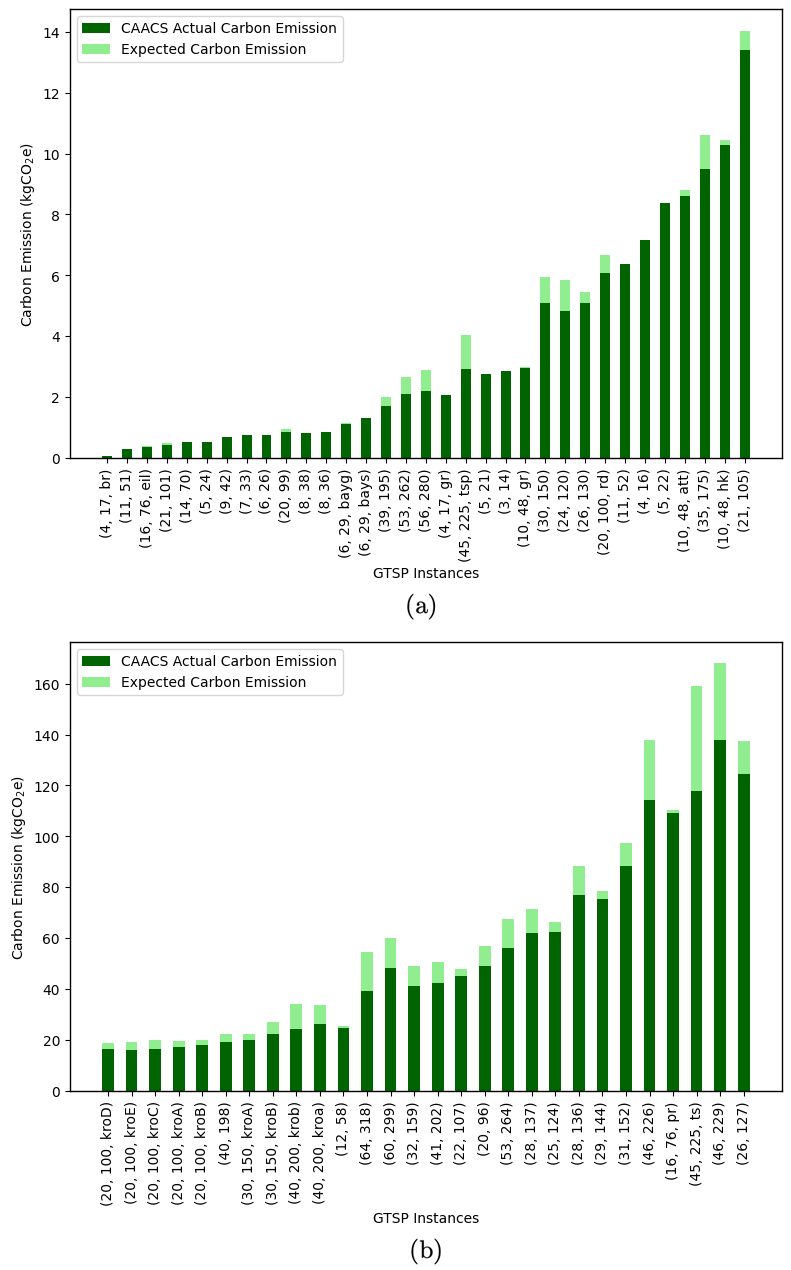}
    \caption{Comparison of expected carbon emission (in light green) in comparison with actual carbon emission. }
    \label{fig:exp_carbon}
\end{figure}

Smaller instances usually result in the same carbon emission,  indicating the same path minimizes the cost and carbon emission as shown in Figure \ref{fig:exp_carbon}. In total, 14 out of the 62 cases had no change in carbon emission. However, 1 case out of the 62 found a carbon emission increase in carbon namely the instance (5, 22) which we attribute to the algorithm's path selection where the optimal path for cost did not align with the optimal path for carbon emission, due to the specific constraints and structure of that instance. Finally, all remaining cases decreased in carbon (47 cases). When comparing the cost and carbon graphs, we find that cases with a tradeoff between cost and carbon meaning the cases in Figure \ref{fig:exp_cost} with a higher error in cost, result in a greater decrease in carbon in Figure \ref{fig:exp_carbon}.  

\begin{figure}[h]
    \centering 
    \hspace*{-3.0in}
    \includegraphics[width=0.42\textwidth]{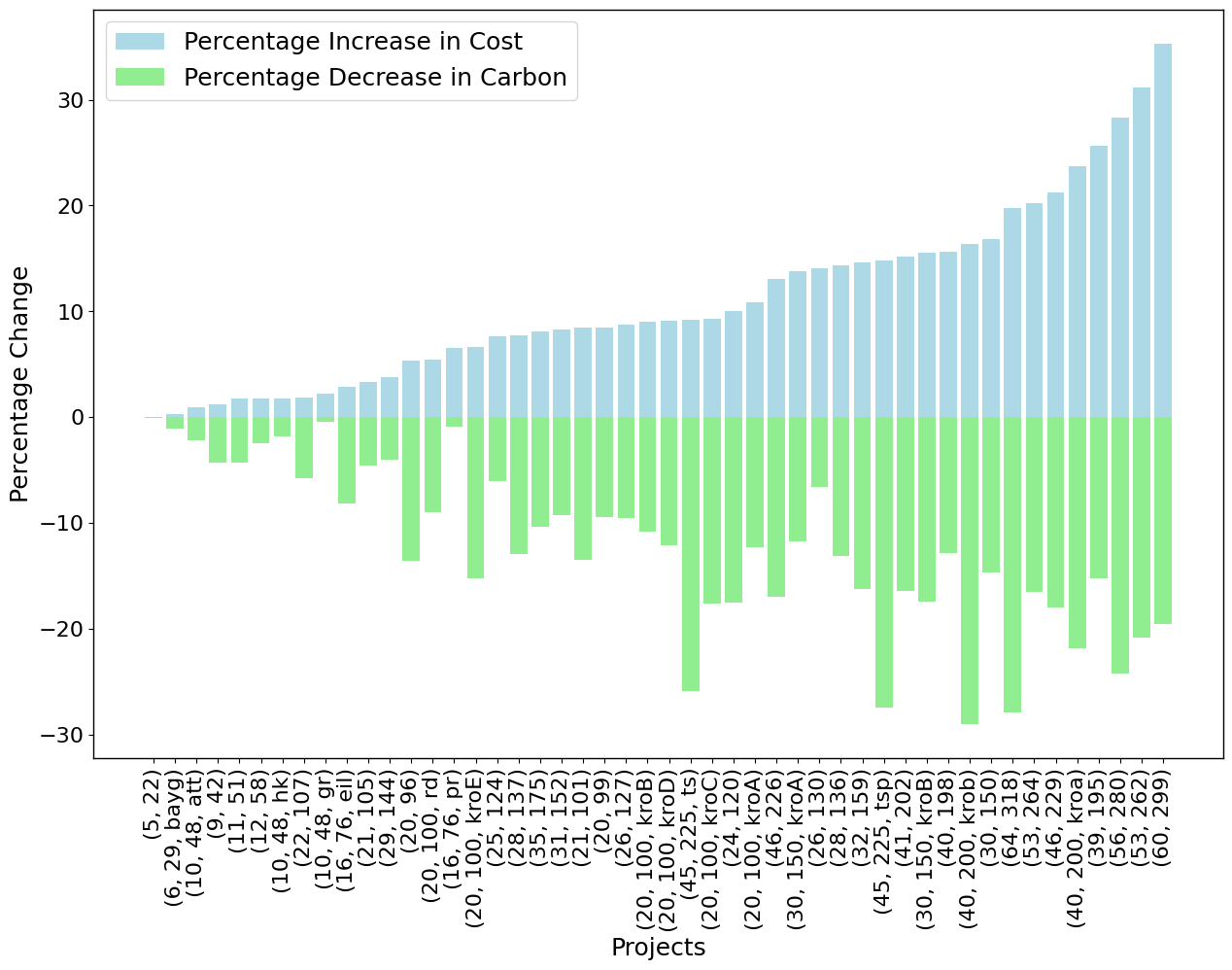}
    \caption{Comparison of the Percentage Increase in Cost and Percentage Decrease in Carbon for Select Project Instances. }
    \label{fig:tradeoff}
\end{figure} 

\medskip Note that a portion of the light green in Figure \ref{fig:exp_carbon} lies under the dark green so the Expected Carbon Emission can be seen as the sum of the dark green and light green. Also, the Expected Emission represents the amount of emission produced corresponding to the optimal cost path. As a result, the light green displayed would represent the amount of carbon emission decreased based on the new path the \textit{CAACS Algorithm} finds. The GTSP instances are ordered in the same order as Figure \ref{fig:exp_cost} for direct comparison. 

\medskip

As the cost increases, there is a corresponding significant decrease in carbon emissions, and this is shown in  Figure \ref{fig:tradeoff}, where one can see there is a strong inverse relationship between the percentage increase in cost and the percentage decrease in carbon emissions.  This suggests that paths that have larger costs and will also reduce carbon emissions (although this doesn't apply for all cases as seen above), are supported by  \cite{DAS2023101816} and \cite{MICHELI2018316}. 

Excluded from Figure \ref{fig:tradeoff}, there were 6 out of the 62 cases with no change in carbon or cost and only 2 cases where there was an increase in cost and no change in carbon.

\section{Time Complexity}
\label{sec: time complexity}

In what follows we shall analyze the time complexity of the \textit{Carbon Aware Ant Colony System Algorithm} by considering the number of nodes $N$ and the number of clusters $M$  first independently,  and then varying $N$ at a fixed ratio to $M$. We measure the number of iterations that the algorithm takes to meet a terminating condition. This condition makes a discovered solution the final solution if it doesn't change after $\frac{N}{5}$ iterations or 1000 iterations are completed. In all trials below, 1000 iterations are never reached. One should note that each iteration of the \textit{CAACS Algorithm} may not follow a linear pattern. However, the exact nature of this non-linearity depends heavily on the specific implementation details, which are outside the scope of this paper. We used $30$ ants, $A = 50$, and $r_0 = 0.5$ to balance the \textit{exploration} and \textit{exploitation phases.} 

In order to understand the time complexity in terms of $N$, let  $M=20$ be fixed. For every   $40\leq N\leq 200$, we generate 10 random graphs by considering the weight matrix instantiated to random values between 10 and 5000, and cluster sizes randomly generated such that the GTSP criteria were met. We ran 10 trials on each of the graphs, for a total of 100 trials for every value of $N$, which were all successful.

\begin{figure}[h]
    \centering 
    \hspace{-2.7in}
    \includegraphics[width=0.35\textwidth]{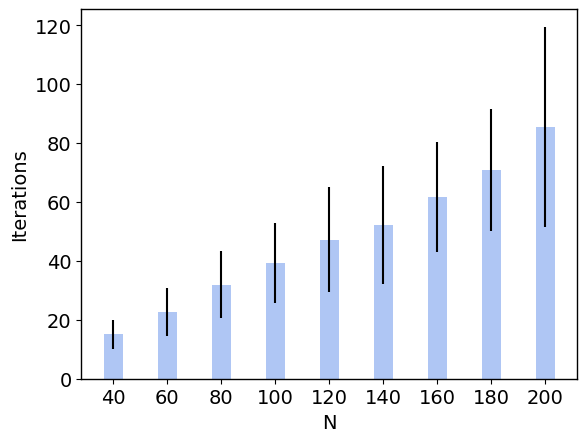}
    \caption{Iterations varying with $N$. Black lines on bars represent standard deviation.}
    \label{fig:fixed_M_tc}
\end{figure}

In Figure \ref{fig:fixed_M_tc}, we see that as $N$ increases the number of iterations appears to increase. This increase appears to be linear, as shown by the consistent upward trend of the bars. Additionally, the error bars indicate that as $N$ increases, the variability in the number of iterations also increases. This could be attributed to the larger number of path possibilities within each cluster and increases the complexity of the problem. 

\begin{figure}[h]
    \centering 
    \hspace{-2.7in}
    \includegraphics[width=0.35\textwidth]{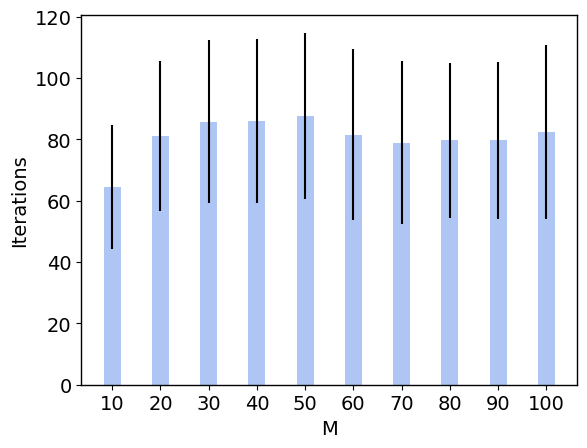}
    \caption{Iterations varying with $M$ with a fixed $N$. Black lines on bars represent standard deviation.}
    \label{fig:fixed_N_tc}
\end{figure}

We shall next analyze the time complexity in terms of $M$, the number of clusters.  Consider $N=200$ fixed and  $10 \leq M\leq 100$. As before, we run 10 trials on 10 graphs for every value of $M$ and use the same termination criteria, leading to Figure \ref{fig:fixed_N_tc}.

The number of iterations appears to be roughly steady after $M = 20$, see Figure \ref{fig:fixed_N_tc}. For $10 \leq M\leq100$, there seems to be a slight increase which can be modeled linearly and the slope of this trend (less than 20 iterations for a 90 unit increase in $M$), in addition to the relatively constant number of iterations as clusters increase, suggests that this algorithm scales well to a larger number of clusters.

Finally, we shall consider what happens when $N$ varies with $M$, setting $N=  M/5$. We ran 10 trials on 10 graphs where   $50\leq N \leq350$ and $M$ are computed according to Eq.~\eqref{eq:eq8}, and all trials were successful using the same termination criteria as above.

\begin{figure}[h]
    \centering 
    \hspace{-2.7in}
    \includegraphics[width=0.35\textwidth]{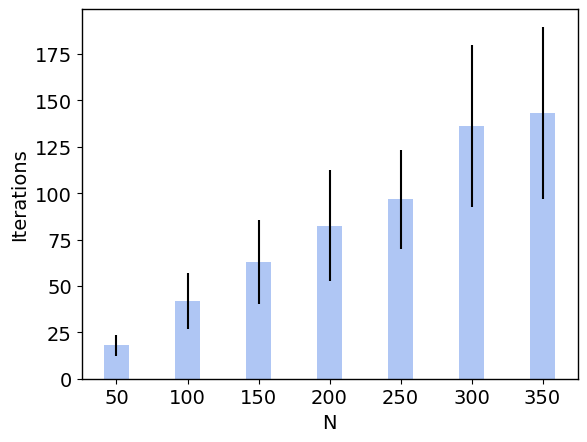}
    \caption{Iterations varying with $N=M/5$, where the blue bar represents the average number of iterations and black lines on bars represent standard deviation.}
    \label{fig:fixed_ratio_tc}
\end{figure}

The number of iterations increases linearly with $N$ as shown in Figure \ref{fig:fixed_ratio_tc}. Similarly, the error increases as the number of nodes increases since there are more possible paths. At $N = 350$, we see the number of iterations is less than the expected number of iterations following a linear regression indicating that the algorithm could scale well to larger problems with more nodes. Increased variability could be due to the random generation of grids resulting in grids with points that are concentrated in out-of-the-way locations. One could argue that at a fixed ratio, the \textit{CAACS Algorithm} is less than linear because the distribution of points ensures that there is sufficient exploration within and between clusters. This balanced exploration helps ants avoid getting stuck in local optima. As a result, the ratio in clustering allows for efficient local and global optimization, reducing the search space, and leading to fewer iterations required even as $N$ increases.

In summary, $M$ can be viewed as a measure of the solution's granularity or detail. As $N$ grows, the run time increases linearly, while increasing the number of clusters $M$ incurs a constant factor. At a fixed ratio, we observed sub-linear time complexity for larger values of $N$. Consequently, with its linear time complexity, we believe the algorithm will scale effectively to handle large problems.

\section{Applications}
\label{sec: applications}

As mentioned before, the GTSP can be used to model a variety of problems among which are microchip design,   DNA sequencing \cite{DNA1, Dundar2019}, delivery and aircraft logistics, network design,  and image retrieval \cite{POP2024819}. To illustrate the novelty of the \textit{CAACS Algorithm} as well as the benefits of its use, we will consider the following:

\begin{itemize}
    \item \textit{\textbf{Network Design.}} We use the \textit{CAACS Algorithm} to develop a sustainable road network in the United States in Section \ref{subsec:network}.
    \item \textit{\textbf{Delivery Logistics.}} We use \textit{CAACS Algorithm} to optimize delivery logistics by finding alternative carbon-aware delivery routes and finding prominent road networks in Section \ref{subsec:delivery}.
    \item \textit{\textbf{Commercial Aircraft Logistics.}} We discuss the application of \textit{CAACS Algorithm} to commercial airplane logistics such as connecting flights and aircraft selection in Section \ref{subsec:planes}.

\end{itemize}

\subsection{Network Design}
\label{subsec:network}

\begin{figure}[h]
    \centering 
    \hspace{-3.7in}
    \includegraphics[width=0.47\textwidth]{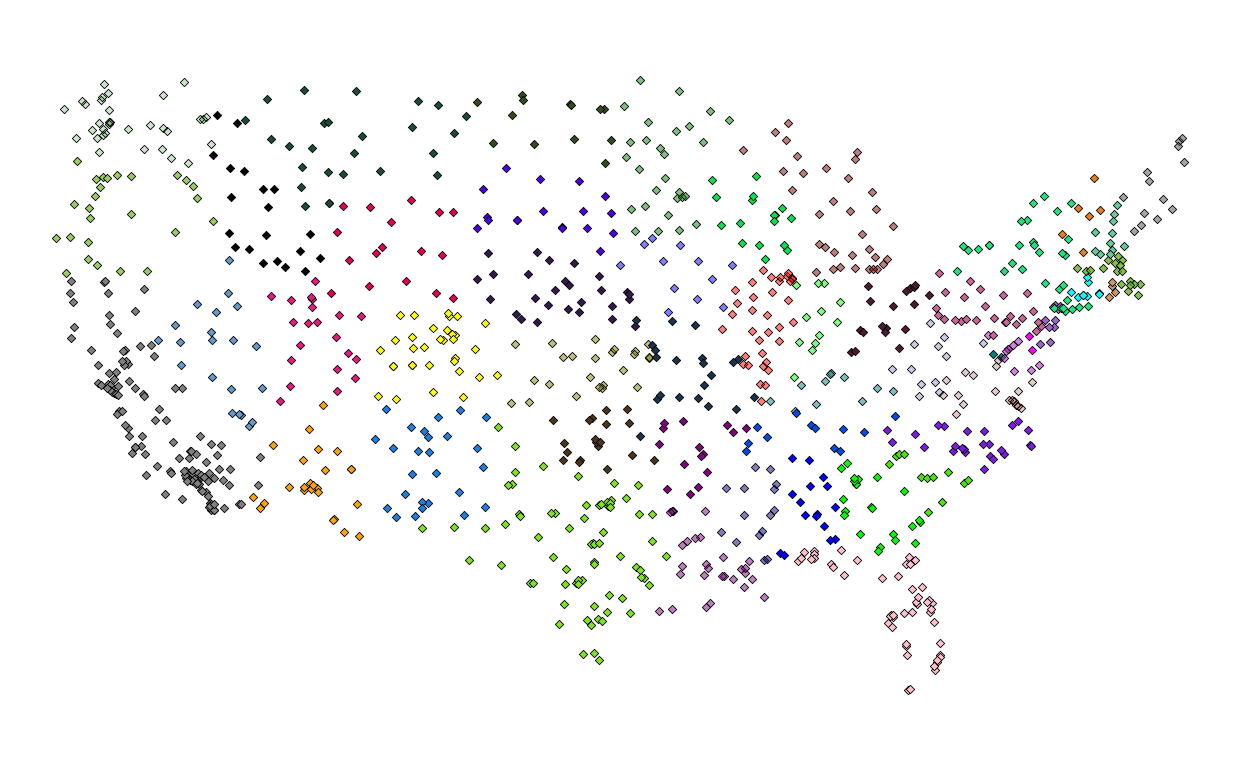}
    \caption{Empty grid with 1097 cities in 48 continental U.S. states (each state a cluster), and DC. }
    \label{fig:usa_empty}
\end{figure}

The \textit{Carbon Aware Ant Colony System Algorithm} can be used to build a road network between   48 continental United States cities  (excluding Alaska and Hawaii) and the District of Columbia (DC). We use the longitude and latitude of the 1097 cities  \cite{helsgaun_glkh}  to generate a rectangular grid of US cities shown in Figure \ref{fig:usa_empty}.

\bigskip

By spawning 20 ants randomly on the graph in Figure \ref{fig:usa_empty},  we ran 10 trails each with 250 iterations, and the final road network is shown below in Figure \ref{fig:usa_final}. 

\bigskip

Note that the paths between nodes in GTSP are straight lines in a theoretical model but when applied to roads, these paths can be translated into actual highway routes. The GTSP algorithm will determine the optimal sequence of nodes, and then this sequence can be mapped onto the highway network, following the roads' actual layout.

\begin{figure}[h]
    \centering 
    \hspace{-3.7in}
    \includegraphics[width=0.47\textwidth]{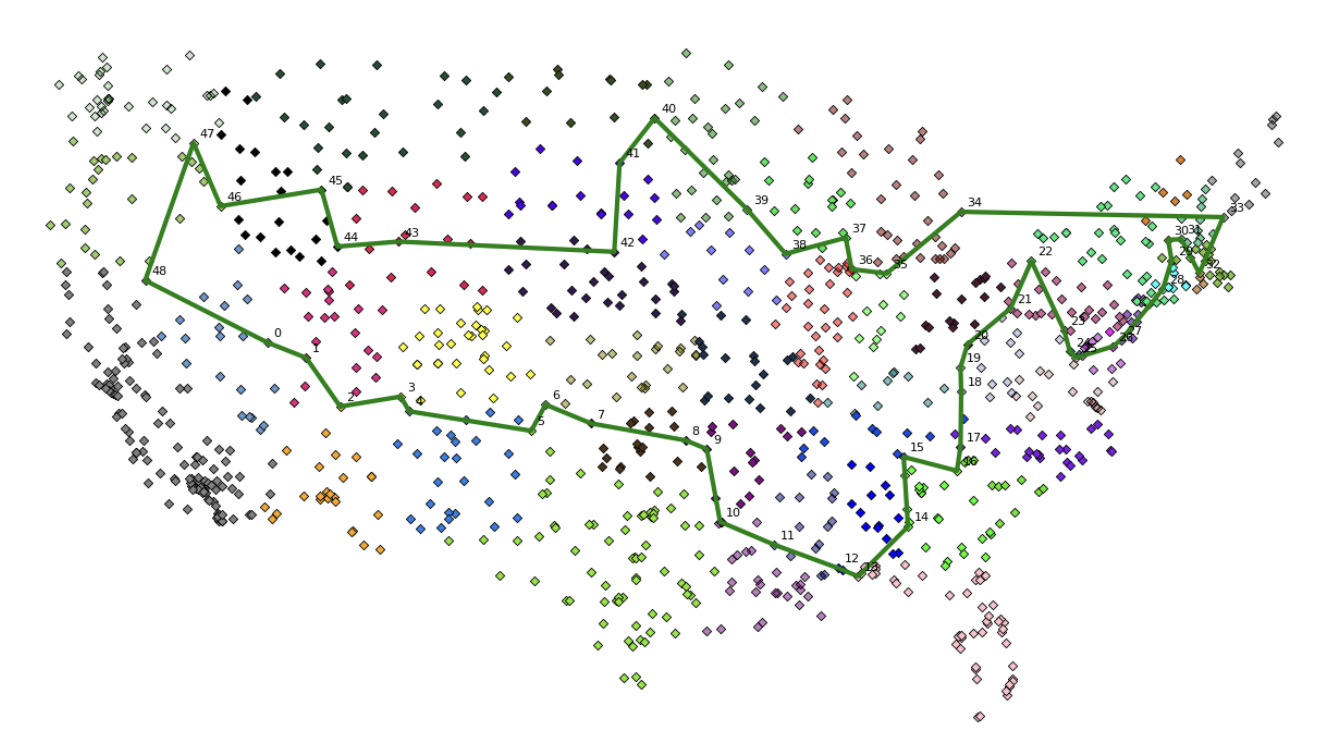}
    \caption{The final road network generated by the algorithm.}
    \label{fig:usa_final}
\end{figure}

\begin{figure}[h]
\begin{center}
    \hspace{-2.5in}
    \includegraphics[width=0.31\textwidth]{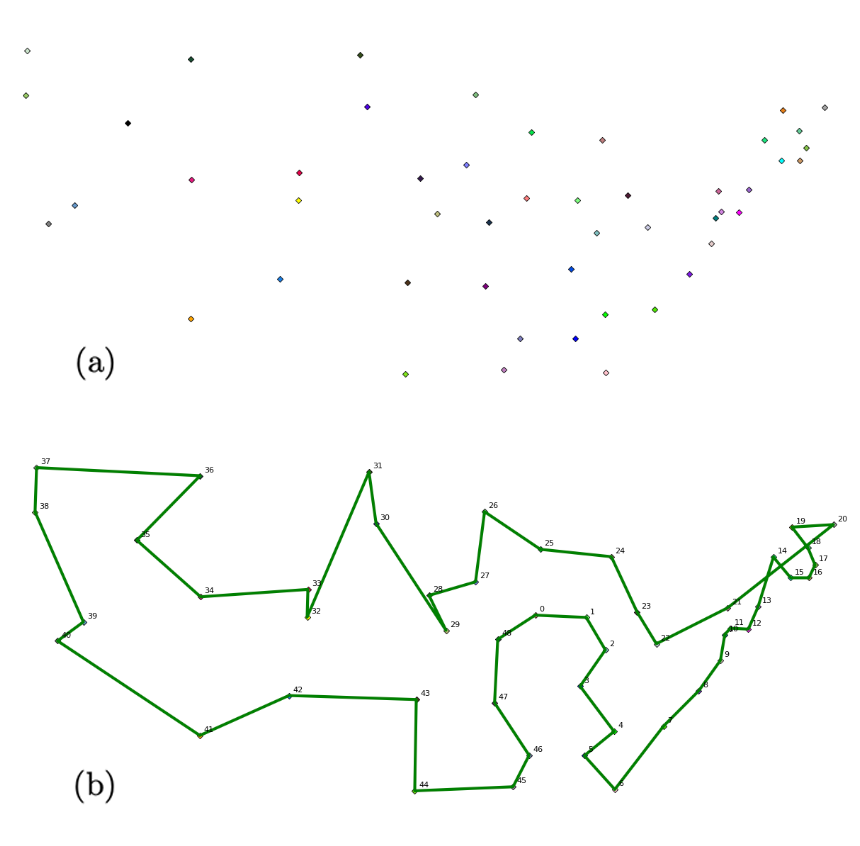}
    \caption{(a). Empty grid with capitals of 49 continental U.S. states (each state a cluster), and DC. (b). The final best network was discovered after 20 trials.}
    \label{fig:og_capital}
    \end{center}
\end{figure}

One can reduce the GTSP to a TSP instance by selecting the capital of every state and finding the least cost path for without carbon and with carbon as seen in Figure \ref{fig:og_capital} (a). Again, we instantiate 20 ants on the graph, and for each trial we ran 500 iterations, running 20 trials in total and selecting the best cost path in Figure \ref{fig:og_capital} (b).

\begin{figure}[h]
  \hspace*{-8.5cm}
    \includegraphics[width=0.48\textwidth]{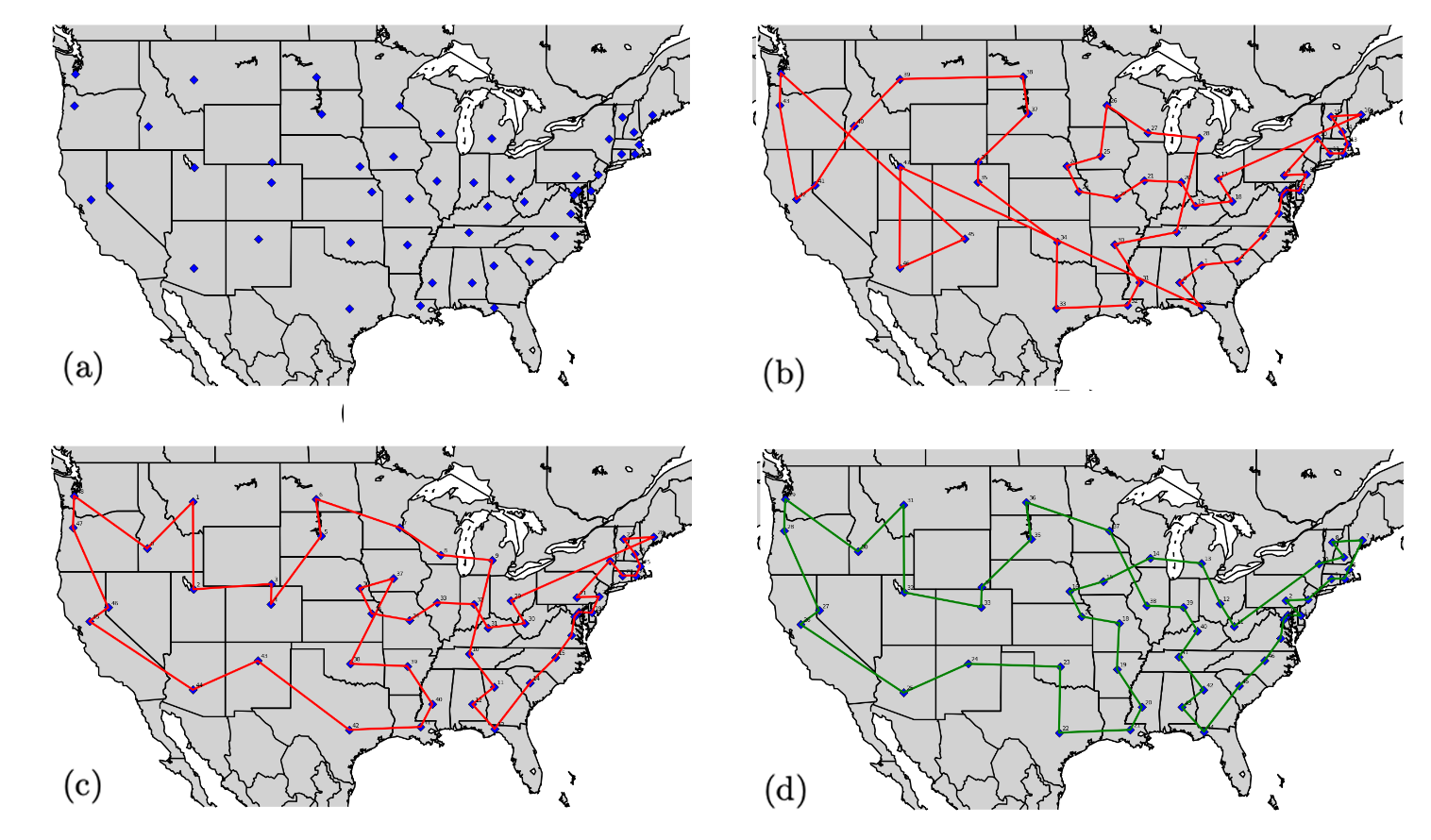}
    \caption{Evolution of the Path Discovery USA Capitals: (a). Initialization of the 48 state capitals and Washington D.C. as blue diamonds. (b). Path Found that reaches all the capitals at the first iteration. (c). Path Found that reaches all the capitals at the 250th iteration. (d). Final valid path was found to the TSP problem after 20 trials. Map created using Basemap.}
    \label{fig:usa_capital_ev}
\end{figure}

We are examining transportation networks that enable companies to transport merchandise (or passengers) and ensure that vehicles return to their original starting points, as vehicles need to return to their starting points. This is why our model offers a valuable tool for identifying such closed paths within the transportation network. The algorithm is particularly suited to the problem of designing large-scale transportation systems since before converging to an optimal tour, the ACO algorithm naturally explores multiple routes. In Figure \ref{fig:usa_capital_ev}, the structure outlines the framework of actual highway routes, adhering to the existing road infrastructure \cite{nhfn_map}.

This approach ensures that the optimized routes are realistic and facilitate efficient freight transportation across the highway network. Additionally, in Figure \ref{fig:usa_capital_ev}, we have a network that can be used for sustainable tourism planning in a road trip to visit all 48 capitals. A potential way for more realistic visualization could be using the APIs in \cite{HieverOptimalRoadTrip, GoogleDistanceMatrixAPI}.

Our algorithm can be applied to many similar problems such as managing waste collection services and urban traffic management. Due to the Ant Colony System's ability to dynamically adapt and optimize routes based on real-time data, it is easily applied to manage traffic flow and reduce congestion in urban environments. Moreover, the \textit{CAACS Algorithm} could be applied to microchip design: the GTSP can be used to optimize the routing of wires and the placement of components on a chip. By considering clusters of connections and ensuring each cluster is visited, the algorithm can reduce the total wire length and improve the overall efficiency of the chip design. Additionally, incorporating carbon-aware factors into the algorithm can lead to environmentally sustainable designs, minimizing the energy consumption required for signal transmission, and the algorithm's time complexity allows it to scale well to large problems.

\subsection{Delivery Routes}
\label{subsec:delivery}

Inspired by \cite{BANIASADI2020444}, and since the \textit{CAACS} model prioritizes sustainable pathways, a natural application is for the creation of environmentally friendly delivery logistics. As an example, we consider UPS, the world's largest package delivery company, with over 40,000 UPS Dropbox locations and over 5000 UPS Stores. Using the data from \cite{ups_facilities_locations_2024}, we focus on the state of Virginia considering all the UPS Stores, UPS Alliance Locations, Authorized Shipping Outlet, and UPS Customer Center for a total of 285 nodes. We exclude the number of UPS Drop Boxes since these locations are typically visited on a less consistent basis compared to other facilities like warehouses and distribution centers. This is because UPS Drop Boxes often serve as convenient drop-off points for customers, and their usage patterns can vary widely depending on factors such as location, time of day, and customer behavior. 

In order to cluster the considered UPS, we unsupervised a machine learning algorithm, K-means, to effectively group locations based on their geographic coordinates. Clustering UPS facility locations enables efficient route planning, resource allocation, and effective management by assigning regional managers to high-density areas.

\begin{figure}[H]
    \centering 
    \hspace*{-0.5cm}
    \includegraphics[width=0.5\textwidth]{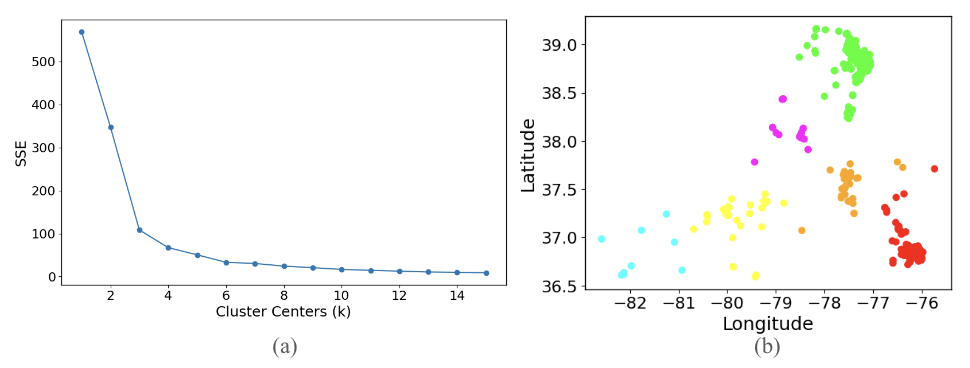}
    \caption{(a). Elbow Plot for Optimal Cluster Selection using K-means, where SSE represents the sum of the squared Euclidean distances of each point to its closest centroid. The red box represents the elbow point region. (b). Resulting 6 Clustered UPS Locations in VA.}
    \label{fig:kmeans}
\end{figure}

\begin{figure}[h]
    \centering 
    \hspace{-2.5in}
    \includegraphics[width=0.38\textwidth]{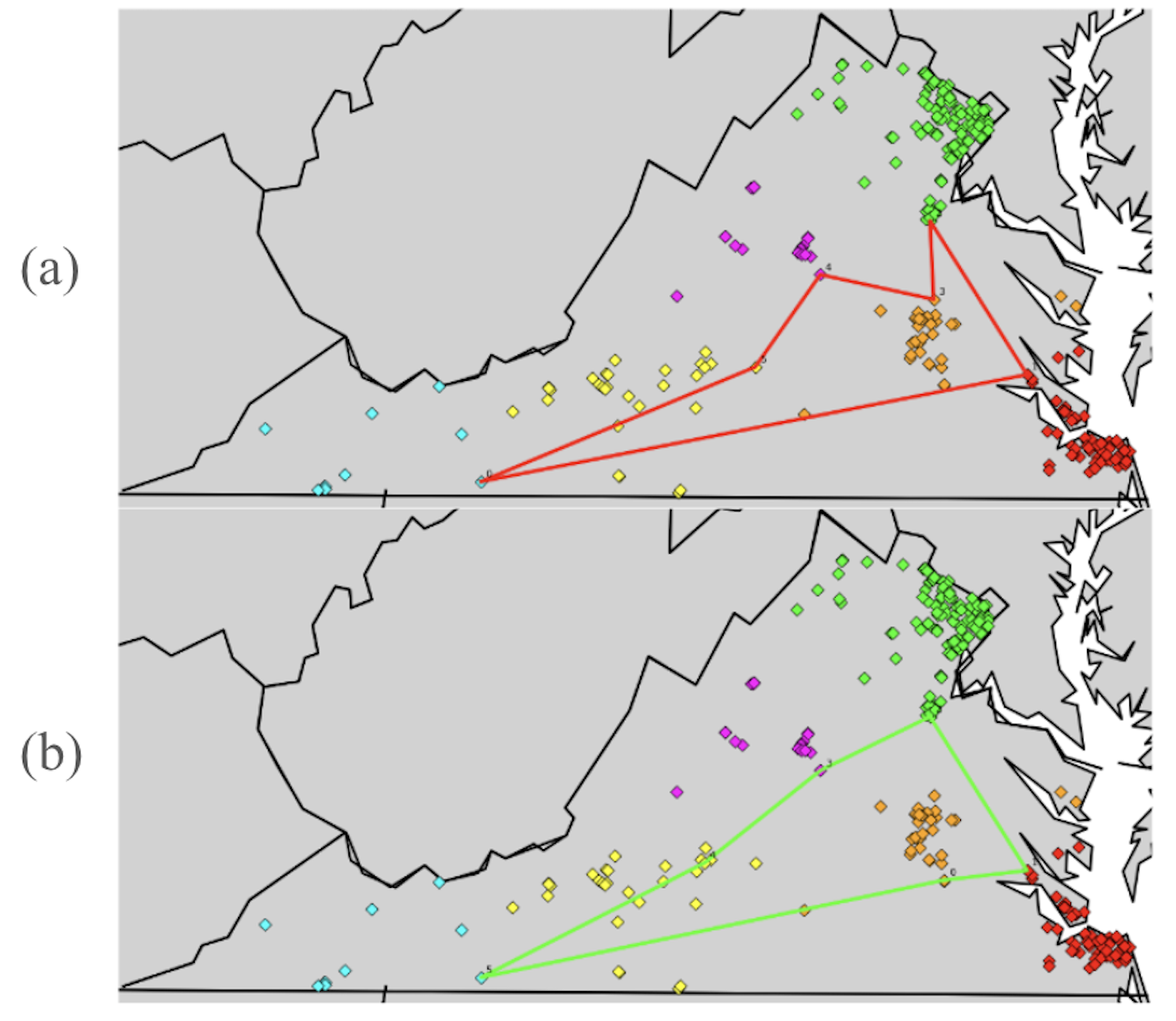}
    \caption{(a). Path discovered after the first iteration of trial with best final path. (b). Best path discovered from all 10 trials. Gray VA Map is created using Baseline. }
    \label{fig:va_6_final}
\end{figure}

As seen in Figure \ref{fig:kmeans} (a), based on the Elbow method below, we select $k = 6$ as the optimal number of clusters. Using this, we ran 10 trials, each with 200 iterations and $r_0 = 0.9$ with the results in Figure \ref{fig:va_6_final}.  This path discovered could be used for more efficient UPS routing. Indeed,  we found a 0.02\% decrease in cost compared to a 1.07\% decrease in Carbon indicating the applicability in finding alternative low-carbon paths in delivery logistics.

\begin{figure}[H]
    \centering 
        \includegraphics[width=0.38\textwidth]{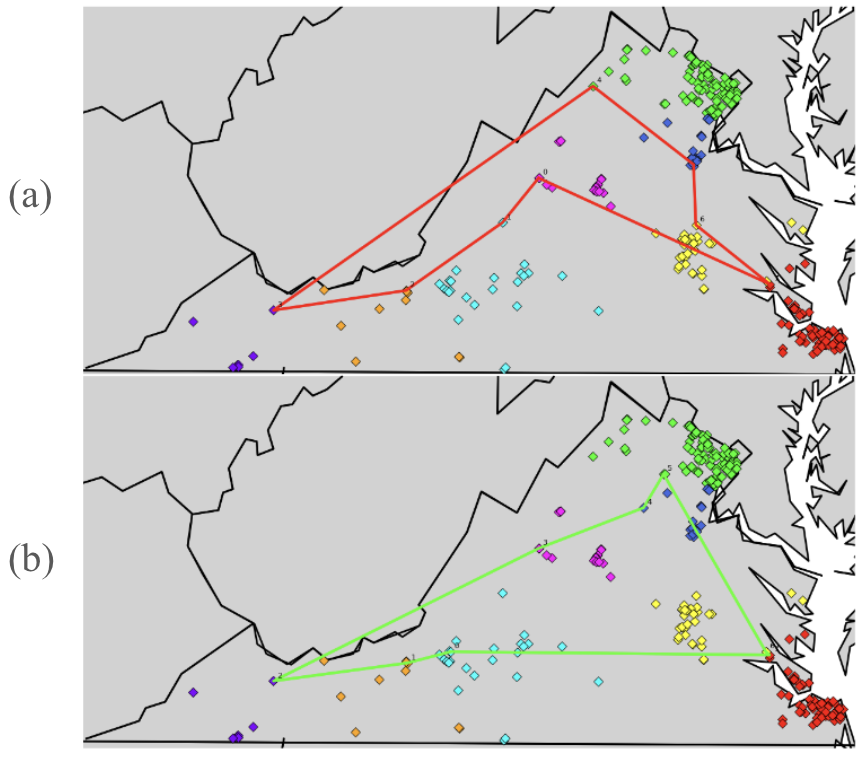}
    \caption{Clustering into 8 groups.}
    \label{fig:8_clusters}
\end{figure}

 This is consistent with our findings in Section \ref{sec:num_ants}, where we find an even larger decrease in carbon than a decrease in cost indicating the promise of alternative routes found by our model. Even though $k = 6$ is the optimal number of clusters according to the elbow method, we add more clusters since some areas have extremely dense clustering (such as Northern VA) indicating a greater population. Because of this, we split the UPS centers into 8 clusters as seen in Figure \ref{fig:8_clusters} and find optimal routes that can visit more population-dense areas. 

\begin{figure}[H]
    \centering 
    \includegraphics[width=0.35\textwidth]{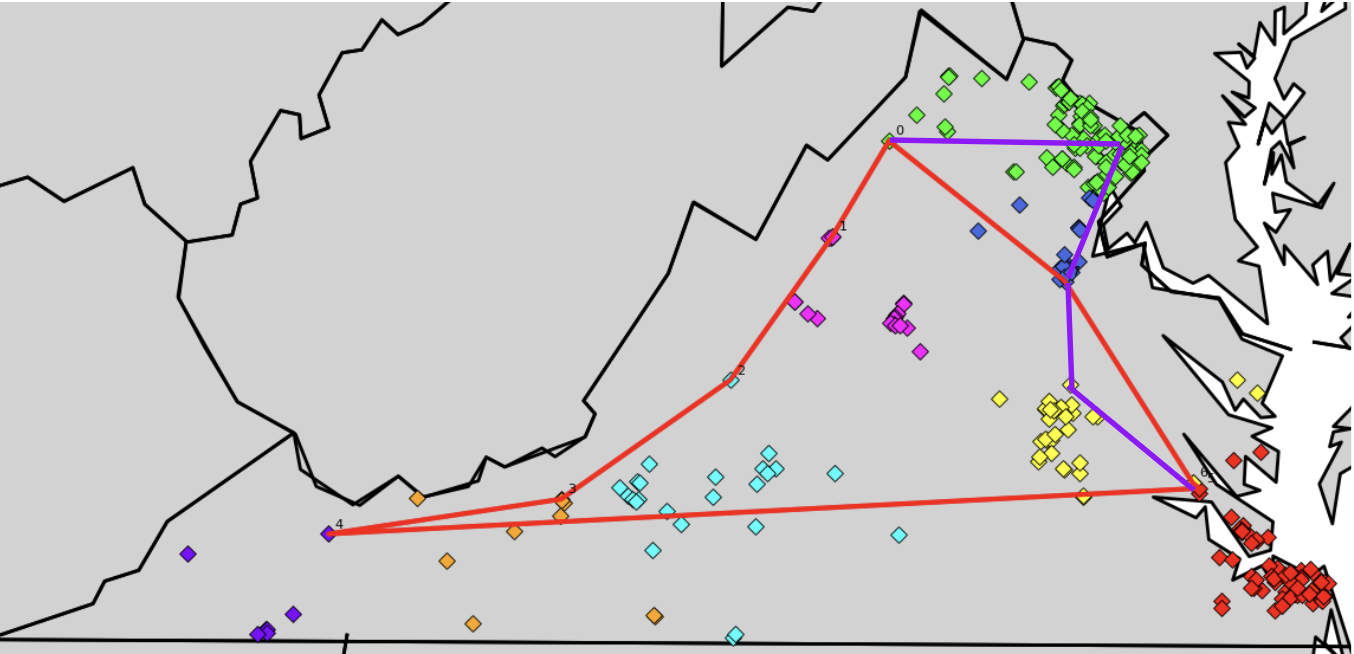}
    \caption{Intermediate route found by the \textit{CAACS Algorithm}. Red Lines represent the path found, and the purple lines represent a slight correction of where the actual interstate is. Different colors represent different clusters.}
    \label{fig:va_8_comparison}
\end{figure}

Additionally, we find parallels between the intermediate paths discovered and the largest Virginia Interstate Highways. For example, in Figure \ref{fig:va_8_comparison}, we see Interstate 81 on the left side being modeled and the correct turning points being selected in the discovered route. The purple lines represent a correction term to where the interstate locations are due to a limited number of clusters in the green sections, where the population is very high. In the yellow cluster, the correction path occurs because the interstate routes pass through the capital, Richmond. This alignment showcases the algorithm's ability to reflect real-world logistics and enhance route optimization in densely populated areas.

The \textit{CAACS Algorithm} can be beneficial for e-commerce companies that need to deliver packages to multiple clusters of destinations (e.g., neighborhoods or regions) with the goal of lowering carbon emissions associated with deliveries. We also see drone-assisted parcel delivery as an application that can be optimized with our algorithm. The algorithm's usage can be easily applied to find sustainable routes for the truck paths to individual depots. Because of the time complexity, it should scale well to a larger problem of delivery routes. 

\subsection{Commercial Aircraft Logistics}
\label{subsec:planes}

Due to the generalizability of this algorithm, it is straightforward to analyze other sources of carbon emissions from different forms of transportation. One of these sectors is the commercial airline industry, which recently entered a carbon-conscious era. To this end, the International Civil Aviation Organization developed the CORSIA program \cite{icao2016carbon}, and the Federal Aviation Administration has also released the United States Aviation Climate Action Plan, which aims to put the aviation industry on a path to net-zero emissions by 2050 \cite{faa_sustainability}. As a result, airlines are forced to consider carbon emission in all parts of their operation process \cite{chao2014assessment}. Additionally, their approaches have broadened to encapsulate multiple factors of both aircraft emission reduction and its economy \cite{hu2022strategies, liu2020drivers}, instead of solely focusing on the economy of aircraft type \cite{liu2017drives, yang2023uncertainty}.

\begin{figure*}[t]
    \flushleft
    \hspace*{-8.9cm}
    \includegraphics[width=1\textwidth]{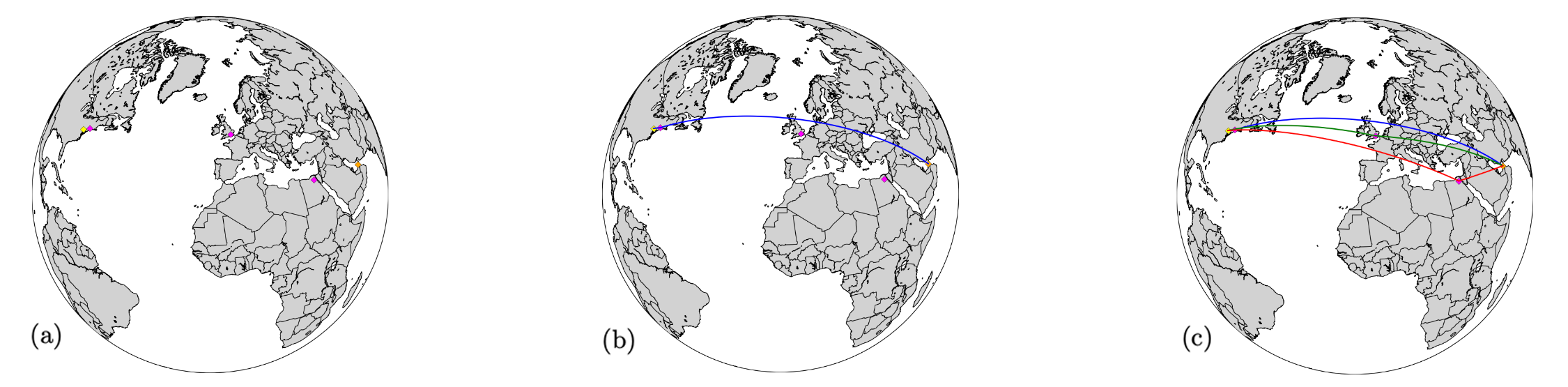}
    \caption{  Aircraft Routes Visualized on the 3D globe (Basemap): (a). Example Travel Locations. Small diamonds represent airports and different clusters are represented using different colors.   (b). Initially discovered path from IAD to DXB and a connection point at EWR. (c). All three potential routes from IAD to DXB with one connecting flight. Different colored paths represent different routes, with the blue path crossing EWR, the green path crossing LHR, and the red path crossing CAI. }
    \label{fig:plane_globe}
\end{figure*}

We shall illustrate how our \textit{CAACS} model can be used in this area by considering the effect of route selection and plane choice on carbon emission by using a case of connecting flights. We shall represent the world as in  Figure \ref{fig:plane_globe} (a), where diamonds represent airports and different colors represent different clusters. In our case, since we are considering an international flight, we only consider international airports; however, this procedure can be easily applied to domestic flights as well.

\begin{itemize}
    \item The yellow diamond represents Dulles International Airport (IAD).
    \item The first pink diamond located in the US represents Newark Liberty International Airport (EWR).
    \item The second pink diamond located in the UK represents Heathrow Airport (LHR). 
    \item The third pink diamond located in Egypt represents Cairo International Airport (CAI). 
    \item The orange diamond represents Dubai International Airport (DXB).
\end{itemize}

Following the GTSP, since only one node is visited per cluster, only one of the three pink diamonds will be visited which accurately models the selection of a connection flight. To generalize the clustering to model connecting flights with $M$ clusters:

\begin{enumerate}
    \item Choose the starting location and set it as Cluster 1 ($C_1$) with only one node. 
    \item Choose the ending location and set it as Cluster M ($C_M$) with only one node.
    \item Find all potential connection flights from $C_1$ and mark all these locations using the same color that is distinct from $C_1 \ldots C_M$. All of these connection flights will be in the same cluster.
    \item Until node in $C_M$ is reached (destination), repeat Step 3, where the next cluster consists of connecting flight locations that can be reached from all nodes in the previous cluster. 
\end{enumerate}

\begin{figure}[h]
    \centering 
    \hspace{-3.3in}
    \includegraphics[width=0.46\textwidth]{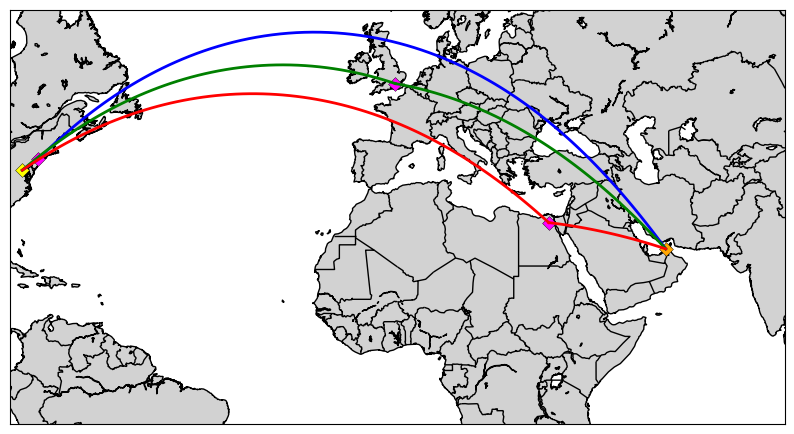}
    \caption{2D Map Representation of the 3 Aircraft Routes in Figure \ref{fig:plane_globe} (c). Map is created using Basemap.}
    \label{fig:route_globe}
\end{figure}

Following the previous procedure,  one can find an optimal one-way route from IAD to DXB, shown in  Figure \ref{fig:plane_globe} and Figure \ref{fig:route_globe}, which represents half of a GTSP solution. However, the GTSP models a full round trip, and the return flight from DXB to IAD can be shown by reversing the starting and ending locations of the trip, finding the new connecting flight locations, and setting conditionals on when a cluster can be visited.

To estimate the carbon emission, we use 2 different models and compare the two results in Table \ref{table:carbon_planes}. The first simple model calculates the distance for each segment of the flight using the Haversine formula, which determines the shortest distance over the Earth's surface between two points. Using this distance, we estimate the carbon emissions using a factor of 0.09 kgCO$_2$ per seat per km \cite{graver2020co2, jardine2009carbon}. The predicted value can be seen in the third column of Table \ref{table:carbon_planes}. This can also be easily updated as planes become more fuel efficient and can also use a more complex carbon function that accounts for different aircraft types, engines, takeoff, weather conditions, etc. 

\begin{table*}[t]
\begin{small}
    \def\arraystretch{1.6}
    \centering
    \begin{tabular}{ | p{3.7cm} | p{2cm} | p{2.5cm} | p{2.5cm} | p{3.3cm} | p{2.8cm} | }
    \hline
    \textbf{Path} & \textbf{Distance (km)} & \textbf{Predicted Carbon} & \textbf{Predicted Average Carbon} & \textbf{Actual Carbon Approx. Range}\\ \hline
   $ IAD \rightarrow EWR \rightarrow DXB $& 11354 km & 1022 kgCO$_2$/seat & 936 kgCO$_2$/seat & 1000--1100 kgCO$_2$/seat\\ \hline
  $  IAD \rightarrow LHR \rightarrow DXB$ & 9360 km & 842 kgCO$_2$/seat & 936 kgCO$_2$/seat & 820--1000 kgCO$_2$/seat \\ \hline
  $  IAD \rightarrow CAI \rightarrow DXB $& 11798 km  & 1062 kgCO$_2$/seat & 936 kgCO$_2$/seat & 810--840 kgCO$_2$/seat \\  \hline
    \end{tabular}
    \caption{Carbon Emission Results for Comparison between Simple Carbon Model and Travel Impact Model. }
    \label{table:carbon_planes}
    \end{small}
\end{table*}

A different carbon function that we shall consider is the Travel Impact Model (TIM) by Google \cite{GoogleTravelImpactModel} that incorporates a large variety of plane types. For each flight, the TIM considers several factors, such as the Great Circle distance between airports and the aircraft type, making it a robust model. Additionally, the TIM estimates fuel burn based on the Tier 3 methodology for emission estimates from the \cite{EEA2019} published by the European Environment Agency and considers:

\begin{itemize}
    \item Tank-to-Wake (TTW) emissions (emissions produced by burning jet fuel during flying, take-off and landing) and 
    \item Well-to-Tank (WTT) emissions (emissions generated during the production, processing, handling, and delivery of jet fuel).
\end{itemize}

In Table \ref{table:carbon_planes}, all carbon values are in (kgCO$_2$/seat) the 4th column consists of the typical carbon emission for this route and the TMI only considers initial and final destinations. The 5th column consists of an approximate range of the actual carbon amounts for a given path from Jul. 4 to Jul 7, a date range that has all 3 types of flights. Although the numbers may not apply to all days, they are shown to illustrate potential causes. A few notable observations are the following:

\begin{enumerate}
    \item Even though IAD $\rightarrow$ CAI $\rightarrow$ DXB (red path in Figure \ref{fig:route_globe}) has the longest distance, it has a lower carbon emissions due to the longer segment from IAD to CAI using the Boeing 787 model which is more carbon efficient. 
    \item  High variability in carbon emission for the IAD $\rightarrow$ LHR $\rightarrow$ DXB (green path in Figure \ref{fig:route_globe}) which is attributed to multiple companies using different types of planes.
    \item In general, United British Airways tended to have a higher range of carbon emissions since they used older Boeing and Airbus Models compared to Virgin Atlantic which used the newer Airbus A350 model. 
\end{enumerate}

After running the \textit{CAACS Algorithm} on this example, we found the final path selected is IAD $\rightarrow$ LHR $\rightarrow$ DXB. After analyzing several days from July 1, 2024, to August 1, 2024, we found out of the three paths IAD $\rightarrow$ LHR $\rightarrow$ DXB had the most flights. Our algorithm is particularly suited to the selection of routes with high flight frequencies, ensuring greater reliability and flexibility for travelers. This path not only maximizes flight availability but also minimizes potential travel disruptions, making it an optimal choice for efficient and dependable air travel. 

\section{Concluding Remarks}
\label{sec:conclusion}

We have dedicated this paper to developing a bi-objective Ant Colony System to solve the GTSP. We have analyzed the effects of the new Scaling Emission Factor and number of ants on the algorithm’s performance before discussing the time complexity and conducted an empirical analysis of our algorithm against benchmark GTSP datasets. Finally, we considered different applications of the \textit{Carbon Aware Ant Colony System Algorithm} to form road networks, to model delivery routing, and to understand sustainability considerations taken by airline companies. \\

 \noindent {\bf Summary of results.} We present a very unique and different approach to finding more sustainable routes to solve the Generalized Traveling Salesman Problem which we hope will inspire further growth and innovation. The current Ant Colony System approach for the TSP was developed in 1997 and since then, there have not been approaches using this ACS model to incorporate environmental sustainability for both the GTSP and the TSP. We hope that the new \textit{Carbon Aware Ant Colony System (CAACS) Algorithm} will spark new ideas in regard to this problem and multi-objective optimization. \\

Additionally, the algorithm offers several particularly interesting advantages. It can be incrementally tuned to find the best path, effectively balancing exploration and exploitation to find a variety of routes. Our algorithm also results in paths that have a higher percentage decrease in carbon emissions compared to a percentage increase in optimal cost which allows it to model real-world routes. This is particularly useful for applications as the algorithm aids in making informed decisions regarding road network creation, delivery logistics, and sustainable travel planning, with a strong emphasis on carbon emissions and cost. Moreover, due to its adaptability to various sources of carbon emissions from different forms of transportation, the algorithm is particularly applicable to commercial aircraft logistics. \\

 \noindent{\bf Future work.} 
The new algorithm developed in this paper could be considered to understand the clustered GTSP, especially considering the application of the problem to drone-assisted parcel delivery \cite{BANIASADI2020444}. Additionally, with the increased projection of drone-assisted delivery usage by commercial companies such as Amazon \cite{drones7030191, cnbc2024amazon}, our algorithm becomes more relevant to real-world simulation. We suggest exploring the 2021 Amazon last mile routing data set \cite{Merchan2022}, which consists of historical routes executed by Amazon delivery drivers. \\
 
Continuing our airplane example in Section \ref{subsec:planes}, we believe the addition of multiple constraints should be added to better model real-world situations. For instance:

\begin{itemize}
    \item \textbf{Time Constraints:} Time is a crucial factor for passengers booking flights. Incorporating departure and arrival time windows, layover durations, and total travel time into the optimization process can ensure that the routes generated by the algorithm are practical for passengers.
    \item \textbf{Social Sustainability:} To ensure social sustainability, the availability and working shift patterns of pilots and airport workers should be considered.
    \item \textbf{Environmental Factors:} Additional environmental factors such as weather conditions and uncertainty should be considered. 
\end{itemize}

 This could perhaps be done by developing heuristics that prioritize routes within acceptable time frames, developing a sub-module within the algorithm to handle compatibility with workforce scheduling, and potentially integrating a weather prediction API to adjust routes based on forecasts. Finally, it would be interesting to see how political factors such as Bilateral Air Service Agreements, that limit airlines from their respective countries can operate, be considered in algorithms \cite{su151612325}. \\

\par On the algorithmic side, future work can be done on trying to optimize our baseline solution. Currently, our models use static update rules that alter the ant graph at the degree until the ending criteria is set which causes solutions to occasionally get stuck at local optima. While local search techniques such as 2-opt or 3-opt could be incorporated within the ACS framework to refine solution quality, it could potentially limit the sustainability aspect. It would be interesting to explore adaptive pheromone updates that dynamically adjust the influence of cost and emissions based on the progress of the algorithm similar to the usage of adaptive learning rates in Machine Learning. \\

\par Finally, from a biological perspective, future work includes examining how different environmental factors affect ant colony behavior and seeing how ants adapt to these stressors. This information could influence the pheromone distributions and be utilized to enhance the \textit{Carbon Aware Ant Colony System Algorithm}. \\

 {\bf Acknowledgments.}  The  authors are thankful to MIT
PRIMES-USA for the opportunity to do this research together and for their invaluable support throughout this research. The work for LPS is partially supported by NSF FRG Award DMS- 2152107, a Simons Fellowship and NSF CAREER Award DMS 1749013. 
 \smallbreak
 {\bf Affiliations.}\\
  (a)  Thomas Jefferson High School for Science and Technology, Alexandria,   USA. \\
  (b)  University of Illinois, Chicago,  USA. \\

\addcontentsline{toc}{section}{%
\protect\numberline{9}%
Bibliography}%
\bibliographystyle{IEEEtran}

\bibliography{IEEEabrv, Primes2024_v3}

\begin{thebibliography}{10}
\providecommand{\url}[1]{#1}
\csname url@samestyle\endcsname
\providecommand{\newblock}{\relax}
\providecommand{\bibinfo}[2]{#2}
\providecommand{\BIBentrySTDinterwordspacing}{\spaceskip=0pt\relax}
\providecommand{\BIBentryALTinterwordstretchfactor}{4}
\providecommand{\BIBentryALTinterwordspacing}{\spaceskip=\fontdimen2\font plus
\BIBentryALTinterwordstretchfactor\fontdimen3\font minus
  \fontdimen4\font\relax}
\providecommand{\BIBforeignlanguage}[2]{{%
\expandafter\ifx\csname l@#1\endcsname\relax
\typeout{** WARNING: IEEEtran.bst: No hyphenation pattern has been}%
\typeout{** loaded for the language `#1'. Using the pattern for}%
\typeout{** the default language instead.}%
\else
\language=\csname l@#1\endcsname
\fi
#2}}
\providecommand{\BIBdecl}{\relax}
\BIBdecl

\bibitem{su9010068}
E.~Eizenberg and Y.~Jabareen, ``Social sustainability: A new conceptual
  framework,'' \emph{Sustainability}, vol.~9, no.~1, 2017.

\bibitem{VALLANCE2011342}
S.~Vallance, H.~C. Perkins, and J.~E. Dixon, ``What is social sustainability? a
  clarification of concepts,'' \emph{Geoforum}, vol.~42, no.~3, pp. 342--348,
  2011, themed Issue: Subaltern Geopolitics.

\bibitem{ELLIOTT2005263}
S.~R. Elliott, ``Sustainability: an economic perspective,'' \emph{Resources,
  Conservation and Recycling}, vol.~44, no.~3, pp. 263--277, 2005,
  sustainability and Renewable Resources.

\bibitem{Ikerd+2012}
J.~Ikerd, \emph{The Essentials of Economic Sustainability}.\hskip 1em plus
  0.5em minus 0.4em\relax Lynne Rienner Publishers, 2012.

\bibitem{environmental1}
J.~C. Little, E.~T. Hester, and C.~C. Carey, ``Assessing and enhancing
  environmental sustainability: A conceptual review,'' \emph{Environmental
  Science $\&$ Technology}, vol.~50, no.~13, pp. 6830--6845, 2016.

\bibitem{environmental2}
A.~Laurent, S.~I. Olsen, and M.~Z. Hauschild, ``Limitations of carbon footprint
  as indicator of environmental sustainability,'' \emph{Environmental Science
  $\&$ Technology}, vol.~46, no.~7, pp. 4100--4108, 2012.

\bibitem{EPA2024}
\BIBentryALTinterwordspacing
{United States Environmental Protection Agency}, ``Carbon pollution from
  transportation,'' 2024. [Online]. Available:
  \url{https://www.epa.gov/transportation-air-pollution-and-climate-change/carbon-pollution-transportation}
\BIBentrySTDinterwordspacing

\bibitem{1976Gt1}
N.~Biggs \emph{et~al.}, \emph{Graph theory 1736-1936}.\hskip 1em plus 0.5em
  minus 0.4em\relax Clarendon, 1976.

\bibitem{Korte2008}
B.~Korte and J.~Vygen, \emph{The Traveling Salesman Problem}.\hskip 1em plus
  0.5em minus 0.4em\relax Berlin, Heidelberg: Springer Berlin Heidelberg, 2008,
  pp. 527--562.

\bibitem{POP2024819}
P.~C. Pop, O.~Cosma, C.~Sabo, and C.~P. Sitar, ``A comprehensive survey on the
  generalized traveling salesman problem,'' \emph{European Journal of
  Operational Research}, vol. 314, no.~3, pp. 819--835, 2024.

\bibitem{KAABACHI2017886}
I.~Kaabachi, D.~Jriji, F.~Madany, and S.~Krichen, ``A bi-criteria ant colony
  optimization for minimizing fuel consumption and cost of the traveling
  salesman problem with time windows,'' \emph{Procedia Computer Science}, vol.
  112, pp. 886--895, 2017, knowledge-Based and Intelligent Information $\&$
  Engineering Systems: Proceedings of the 21st International Conference,
  KES-20176-8 September 2017, Marseille, France.

\bibitem{DAS2023101816}
M.~Das, A.~Roy, S.~Maity, and S.~Kar, ``A quantum-inspired ant colony
  optimization for solving a sustainable four-dimensional traveling salesman
  problem under type-2 fuzzy variable,'' \emph{Advanced Engineering
  Informatics}, vol.~55, p. 101816, 2023.

\bibitem{MICHELI2018316}
G.~J. Micheli and F.~Mantella, ``Modelling an environmentally-extended
  inventory routing problem with demand uncertainty and a heterogeneous fleet
  under carbon control policies,'' \emph{International Journal of Production
  Economics}, vol. 204, pp. 316--327, 2018.

\bibitem{DNA1}
M.~Caserta and S.~Voß, ``A hybrid algorithm for the dna sequencing problem,''
  \emph{Discrete Applied Mathematics}, vol. 163, pp. 87--99, 2014,
  matheuristics 2010.

\bibitem{Dundar2019}
A.~O. Dundar, M.~A. Sahman, M.~Tekin, and M.~S. Kıran, ``A comparative
  application regarding the effects of traveling salesman problem on logistics
  costs,'' \emph{International Journal of Intelligent Systems and Applications
  in Engineering}, vol.~7, no.~4, pp. 207--2015, Dec. 2019.

\bibitem{LIEN1993177}
Y.-N. Lien, E.~Ma, and B.~W.-S. Wah, ``Transformation of the generalized
  traveling-salesman problem into the standard traveling-salesman problem,''
  \emph{Information Sciences}, vol.~74, no.~1, pp. 177--189, 1993.

\bibitem{SMITH20171}
S.~L. Smith and F.~Imeson, ``Glns: An effective large neighborhood search
  heuristic for the generalized traveling salesman problem,'' \emph{Computers
  $\&$ Operations Research}, vol.~87, pp. 1--19, 2017.

\bibitem{Matai10}
R.~Matai, S.~Singh, and M.~L. Mittal, ``Traveling salesman problem: an overview
  of applications, formulations, and solution approaches,'' in \emph{Traveling
  Salesman Problem}.\hskip 1em plus 0.5em minus 0.4em\relax Rijeka: IntechOpen,
  2010, ch.~1.

\bibitem{10.5555/2209505.2209517}
I.~Kara, H.~Guden, and O.~N. Koc, ``New formulations for the generalized
  traveling salesman problem,'' in \emph{New formulations for the generalized
  traveling salesman problem}.\hskip 1em plus 0.5em minus 0.4em\relax World
  Scientific and Engineering Academy and Society (WSEAS), 2012.

\bibitem{pop2007}
P.~Pop, ``New integer programming formulations of the generalized travelling
  salesman problem,'' \emph{American Journal of Applied Sciences}, vol.~4, pp.
  932--937, 11 2007.

\bibitem{henry1969record}
A.~Henry-Labordere, ``The record balancing problem: A dynamic programming
  solution of a generalized traveling salesman problem rairo, vol,'' \emph{The
  record balancing problem: A dynamic programming solution of a generalized
  traveling salesman problem RAIRO vol}, 1969.

\bibitem{transformation}
\BIBentryALTinterwordspacing
C.~E. Noon and J.~C. Bean, ``An efficient transformation of the generalized
  traveling salesman problem,'' \emph{INFOR: Information Systems and
  Operational Research}, vol.~31, no.~1, pp. 39--44, 1993. [Online]. Available:
  \url{https://doi.org/10.1080/03155986.1993.11732212}
\BIBentrySTDinterwordspacing

\bibitem{reduction}
G.~Gutin and D.~Karapetyan, ``Generalized traveling salesman problem reduction
  algorithms,'' \emph{Algorithmic Operations Research}, vol.~4, no.~2, pp.
  144--154, 2009.

\bibitem{approx}
P.~C. Pop, O.~Cosma, C.~Sabo, and C.~P. Sitar, ``A comprehensive survey on the
  generalized traveling salesman problem,'' \emph{European Journal of
  Operational Research}, vol. 314, no.~3, pp. 819--835, 2024.

\bibitem{585892}
M.~Dorigo and L.~Gambardella, ``Ant colony system: a cooperative learning
  approach to the traveling salesman problem,'' \emph{IEEE Transactions on
  Evolutionary Computation}, vol.~1, no.~1, pp. 53--66, 1997.

\bibitem{pintea2017generalized}
C.-M. Pintea, P.~C. Pop, and C.~Chira, ``The generalized traveling salesman
  problem solved with ant algorithms,'' \emph{Complex Adaptive Systems
  Modeling}, vol.~5, no.~1, p.~8, 2017.

\bibitem{Wang2018}
D.~Wang, D.~Tan, and L.~Liu, ``Particle swarm optimization algorithm: an
  overview,'' \emph{Soft Computing}, vol.~22, no.~2, pp. 387--408, Jan 2018.

\bibitem{ACO}
M.~Dorigo, M.~Birattari, and T.~Stutzle, ``Ant colony optimization,''
  \emph{IEEE Computational Intelligence Magazine}, vol.~1, no.~4, pp. 28--39,
  2006.

\bibitem{dorigo}
A.~C. Marco~Dorigo, Vittorio~Maniezzo, ``Ant system: optimization by a colony
  of cooperating agents,'' \emph{IEEE Transactions on Systems, Man, and
  Cybernetics, Part B (Cybernetics)}, vol.~26, pp. 29--41, 1996.

\bibitem{czaczkes2011synergy}
T.~J. Czaczkes, C.~Grüter, S.~M. Jones, and F.~L.~W. Ratnieks, ``Synergy
  between social and private information increases foraging efficiency in
  ants,'' \emph{Biology Letters}, vol.~7, no.~4, pp. 521--524, 2011.

\bibitem{ant_col_video1}
H.~Lam, ``Ant colony optimization - part 3.1: Simple ant colony optimization
  (saco),'' Youtube, 2022.

\bibitem{ant_system}
M.~Dorigo, V.~Maniezzo, and A.~Colorni, ``Ant system: optimization by a colony
  of cooperating agents,'' \emph{IEEE Transactions on Systems, Man, and
  Cybernetics, Part B (Cybernetics)}, vol.~26, no.~1, pp. 29--41, 1996.

\bibitem{JUNMAN2012319}
K.~Jun-man and Z.~Yi, ``Application of an improved ant colony optimization on
  generalized traveling salesman problem,'' \emph{Energy Procedia}, vol.~17,
  pp. 319--325, 2012, 2012 International Conference on Future Electrical Power
  and Energy System.

\bibitem{app132111817}
B.~Ghimire, A.~Mahmood, and K.~Elleithy, ``Hybrid parallel ant colony
  optimization for application to quantum computing to solve large-scale
  combinatorial optimization problems,'' \emph{Applied Sciences}, vol.~13,
  no.~21, 2023.

\bibitem{Karapetyan2012AnEH}
D.~Karapetyan and M.~Reihaneh, ``An efficient hybrid ant colony system for the
  generalized traveling salesman problem,'' \emph{Algorithmic Oper. Res.},
  vol.~7, pp. 22--29, 2012.

\bibitem{Pintea2017}
C.-M. Pintea, P.~C. Pop, and C.~Chira, ``The generalized traveling salesman
  problem solved with ant algorithms,'' \emph{Complex Adaptive Systems
  Modeling}, vol.~5, no.~1, p.~8, Aug 2017.

\bibitem{su152115457}
\BIBentryALTinterwordspacing
F.~Jelti, A.~Allouhi, and K.~A. Tabet~Aoul, ``Transition paths towards a
  sustainable transportation system: A literature review,''
  \emph{Sustainability}, vol.~15, no.~21, 2023. [Online]. Available:
  \url{https://www.mdpi.com/2071-1050/15/21/15457}
\BIBentrySTDinterwordspacing

\bibitem{ICCT2018}
\BIBentryALTinterwordspacing
I.~C. on~Clean~Transportation, ``Beyond road zevs: Considerations for
  zero-emission vehicles in the off-road sector,'' International Council on
  Clean Transportation, Tech. Rep., 2018. [Online]. Available:
  \url{https://theicct.org/sites/default/files/publications/Beyond_Road_ZEV_Working_Paper_20180718.pdf}
\BIBentrySTDinterwordspacing

\bibitem{EDF2021}
\BIBentryALTinterwordspacing
E.~D. Fund, ``Feasibility of electrifying medium- and heavy-duty vehicles,''
  Environmental Defense Fund, Tech. Rep., 2021. [Online]. Available:
  \url{https://www.edf.org/sites/default/files/documents/EDFMHDVEVFeasibilityReport22jul21.pdf}
\BIBentrySTDinterwordspacing

\bibitem{NHTSA2015}
\BIBentryALTinterwordspacing
N.~H. T.~S. Administration, ``National traffic speeds survey iii,'' National
  Highway Traffic Safety Administration, Tech. Rep., 2015. [Online]. Available:
  \url{https://www.nhtsa.gov/sites/nhtsa.gov/files/documents/812485_national-traffic-speeds-survey-iii-2015.pdf}
\BIBentrySTDinterwordspacing

\bibitem{FHWA2013}
\BIBentryALTinterwordspacing
{Federal Highway Administration}, ``Setting speed limits for safety,''
  \emph{Public Roads}, vol.~77, no.~2, 2013. [Online]. Available:
  \url{https://highways.dot.gov/public-roads/septemberoctober-2013/setting-speed-limits-safety}
\BIBentrySTDinterwordspacing

\bibitem{neighborhoodroadside2024}
\BIBentryALTinterwordspacing
{Neighborhood Roadside}, ``What is payload capacity?'' 2024. [Online].
  Available: \url{https://neighborhoodroadside.com/what-is-payload-capacity/}
\BIBentrySTDinterwordspacing

\bibitem{cleantechnol3020028}
\BIBentryALTinterwordspacing
C.~Cunanan, M.-K. Tran, Y.~Lee, S.~Kwok, V.~Leung, and M.~Fowler, ``A review of
  heavy-duty vehicle powertrain technologies: Diesel engine vehicles, battery
  electric vehicles, and hydrogen fuel cell electric vehicles,'' \emph{Clean
  Technologies}, vol.~3, no.~2, pp. 474--489, 2021. [Online]. Available:
  \url{https://www.mdpi.com/2571-8797/3/2/28}
\BIBentrySTDinterwordspacing

\bibitem{engelbrecht2007computational}
A.~P. Engelbrecht, \emph{Computational intelligence: an introduction}.\hskip
  1em plus 0.5em minus 0.4em\relax John Wiley \& Sons, 2007.

\bibitem{constantinou2010ant}
D.~Constantinou \emph{et~al.}, ``Ant colony optimisation algorithms for solving
  multi-objective power-aware metrics for mobile ad hoc networks,'' Ph.D.
  dissertation, University of Pretoria, 2010.

\bibitem{gutin2009memetic}
G.~Gutin and D.~Karapetyan, ``A memetic algorithm for the generalized traveling
  salesman problem,'' \emph{Natural Computing}, vol.~9, no.~1, pp. 47--60,
  2009.

\bibitem{helsgaun_glkh}
K.~Helsgaun, ``Glkh: Generalized large-scale kernel-based heuristic,''
  \url{http://akira.ruc.dk/~keld/research/GLKH/}.

\bibitem{nhfn_map}
{Federal Highway Administration}, ``National highway freight network map,''
  \url{https://ops.fhwa.dot.gov/freight/infrastructure/nfn/maps/nhfn_map.htm},
  2024.

\bibitem{HieverOptimalRoadTrip}
\BIBentryALTinterwordspacing
R.~S. Olson, ``Computing the optimal road trip across the u.s.'' 2015.
  [Online]. Available:
  \url{https://github.com/rhiever/Data-Analysis-and-Machine-Learning-Projects/blob/master/optimal-road-trip/Computing%20the%20optimal%20road%20trip%20across%20the%20U.S..ipynb}
\BIBentrySTDinterwordspacing

\bibitem{GoogleDistanceMatrixAPI}
\BIBentryALTinterwordspacing
G.~Developers, ``Google maps distance matrix api,'' 2024. [Online]. Available:
  \url{https://developers.google.com/maps/documentation/distance-matrix}
\BIBentrySTDinterwordspacing

\bibitem{BANIASADI2020444}
P.~Baniasadi, M.~Foumani, K.~Smith-Miles, and V.~Ejov, ``A transformation
  technique for the clustered generalized traveling salesman problem with
  applications to logistics,'' \emph{European Journal of Operational Research},
  vol. 285, no.~2, pp. 444--457, 2020.

\bibitem{ups_facilities_locations_2024}
\BIBentryALTinterwordspacing
H.~I. Foundation, ``Ups facilities locations dataset,'' 2024. [Online].
  Available:
  \url{https://www.kaggle.com/datasets/thedevastator/ups-facilities-locations-dataset/data}
\BIBentrySTDinterwordspacing

\bibitem{icao2016carbon}
ICAO, ``Carbon offsetting and reduction scheme for international aviation
  (corsia),'' 2016.

\bibitem{faa_sustainability}
\BIBentryALTinterwordspacing
F.~A. Administration, ``Sustainability,'' 2021, accessed: 2024-06-15. [Online].
  Available:
  \url{https://www.faa.gov/sustainability#:~:text=2021%2C%20U.S.%20Transportation%20Sec.,large%20part%20of%20the%20solution.}
\BIBentrySTDinterwordspacing

\bibitem{chao2014assessment}
C.-C. Chao, ``Assessment of carbon emission costs for air cargo
  transportation,'' \emph{Transportation Research Part D: Transport and
  Environment}, vol.~33, pp. 186--195, 2014.

\bibitem{hu2022strategies}
Y.-J. Hu, L.~Yang, H.~Cui, H.~Wang, C.~Li, and B.-J. Tang, ``Strategies to
  mitigate carbon emissions for sustainable aviation: A critical review from a
  life-cycle perspective,'' \emph{Sustainable Production and Consumption},
  vol.~33, pp. 788--808, 2022.

\bibitem{liu2020drivers}
X.~Liu, Y.~Hang, Q.~Wang, and D.~Zhou, ``Drivers of civil aviation carbon
  emission change: A two-stage efficiency-oriented decomposition approach,''
  \emph{Transportation Research Part D: Transport and Environment}, vol.~89, p.
  102612, 2020.

\bibitem{liu2017drives}
X.~Liu, D.~Zhou, P.~Zhou, and Q.~Wang, ``What drives co2 emissions from
  china’s civil aviation? an exploration using a new generalized pda
  method,'' \emph{Transportation Research Part A: Policy and Practice},
  vol.~99, pp. 30--45, 2017.

\bibitem{yang2023uncertainty}
L.~Yang, Y.-J. Hu, H.~Wang, C.~Li, B.-J. Tang, B.~Wang, and H.~Cui,
  ``Uncertainty quantification of co2 emissions from china's civil aviation
  industry to 2050,'' \emph{Journal of Environmental Management}, vol. 336, p.
  117624, 2023.

\bibitem{graver2020co2}
\BIBentryALTinterwordspacing
B.~Graver, D.~Zhang, and D.~Rutherford, ``Co2 emissions from commercial
  aviation: 2013, 2018, and 2019,'' 2020. [Online]. Available:
  \url{https://theicct.org/wp-content/uploads/2021/06/CO2-commercial-aviation-oct2020.pdf}
\BIBentrySTDinterwordspacing

\bibitem{jardine2009carbon}
\BIBentryALTinterwordspacing
C.~N. Jardine, ``Calculating the carbon dioxide emissions of flights,''
  \emph{Grandstand Central}, 2009. [Online]. Available:
  \url{https://www.grandstandcentral.com/research/energy/downloads/jardine09-carboninflights.pdf}
\BIBentrySTDinterwordspacing

\bibitem{GoogleTravelImpactModel}
\BIBentryALTinterwordspacing
Google, ``Travel impact model,'' 2024, accessed: 2024-06-16. [Online].
  Available: \url{https://github.com/google/travel-impact-model/tree/main}
\BIBentrySTDinterwordspacing

\bibitem{EEA2019}
\BIBentryALTinterwordspacing
{European Environment Agency}, \emph{EMEP/EEA Air Pollutant Emission Inventory
  Guidebook 2019}, 2019, accessed: 2024-06-16. [Online]. Available:
  \url{https://www.eea.europa.eu/publications/emep-eea-guidebook-2019/part-b-sectoral-guidance-chapters/1-energy/1-a-combustion/1-a-3-a-aviation/view}
\BIBentrySTDinterwordspacing

\bibitem{drones7030191}
\BIBentryALTinterwordspacing
X.~Li, J.~Tupayachi, A.~Sharmin, and M.~Martinez~Ferguson, ``Drone-aided
  delivery methods, challenge, and the future: A methodological review,''
  \emph{Drones}, vol.~7, no.~3, 2023. [Online]. Available:
  \url{https://www.mdpi.com/2504-446X/7/3/191}
\BIBentrySTDinterwordspacing

\bibitem{cnbc2024amazon}
\BIBentryALTinterwordspacing
``Amazon gets faa approval for drone delivery in california and texas,''
  \emph{CNBC}, May 2024. [Online]. Available:
  \url{https://www.cnbc.com/2024/05/30/amazon-drone-delivery-faa-approval.html#:~:text=Last%20month%2C%20Amazon%20said%20it,the%20end%20of%20the%20decade}
\BIBentrySTDinterwordspacing

\bibitem{Merchan2022}
\BIBentryALTinterwordspacing
D.~Merchan, J.~Arora, J.~Pachon, K.~Konduri, M.~Winkenbach, S.~Parks, and
  J.~Noszek, ``2021 amazon last mile routing research challenge: Data set,''
  \emph{Transportation Science}, 2022. [Online]. Available:
  \url{https://www.amazon.science/publications/2021-amazon-last-mile-routing-research-challenge-data-set}
\BIBentrySTDinterwordspacing

\bibitem{su151612325}
\BIBentryALTinterwordspacing
R.~J. Raimundo, M.~E. Baltazar, and S.~P. Cruz, ``Sustainability in the
  airports ecosystem: A literature review,'' \emph{Sustainability}, vol.~15,
  no.~16, 2023. [Online]. Available:
  \url{https://www.mdpi.com/2071-1050/15/16/12325}
\BIBentrySTDinterwordspacing

\end{thebibliography}

\end{document}